\documentclass{article}

\usepackage{microtype}
\usepackage{graphicx}
\usepackage{booktabs}

\usepackage{hyperref}

\usepackage{listings}
\usepackage{xcolor}

\usepackage[accepted]{icml2025}

\usepackage{amsmath}

\usepackage{subcaption}
\usepackage{amssymb}
\usepackage{multirow}
\usepackage{mathtools}
\usepackage{amsthm}

\usepackage[capitalize,noabbrev]{cleveref}

\theoremstyle{plain}

\theoremstyle{definition}

\theoremstyle{remark}

\usepackage[textsize=tiny]{todonotes}

\icmltitlerunning{Code Retrieval for MILP Instance Generation}

\begin{document}

\twocolumn[
\icmltitle{Code Retrieval for MILP Instance Generation}




\begin{icmlauthorlist}
\icmlauthor{Tianxing Yang}{yyy,comp}
\icmlauthor{Huigen Ye}{yyy,comp}
\icmlauthor{Hua Xu}{yyy,comp}

\end{icmlauthorlist}

\icmlaffiliation{yyy}{Department of Computer Science and Technology, Tsinghua University}
\icmlaffiliation{comp}{Beijing National Research Center for Information Science and Technology}

\icmlcorrespondingauthor{Hua Xu}{xuhua@mail.tsinghua.edu.cn}

\icmlkeywords{Machine Learning, ICML}

\vskip 0.3in
]



\printAffiliationsAndNotice{}  

\begin{abstract}
Mixed-Integer Linear Programming (MILP) is widely used in fields such as scheduling, logistics, and planning. Enhancing the performance of MILP solvers, particularly learning-based solvers, requires substantial amounts of high-quality data. However, existing methods for MILP instance generation typically necessitate training a separate model for each problem class and are computationally intensive when generating new instances. To address these limitations, we reformulate the \emph{MILP Instance Generation} task as \emph{MILP Code Generation} task, enabling efficient, flexible, and interpretable instance generation through code. Since MILP instances generated from code can vary significantly in scale, we introduce \emph{MILP-EmbedSim}, a new similarity metric that accurately measures the similarity between instances of varying sizes within the same problem class. Leveraging this metric, we propose \emph{MILP-Retrieval}, a pipeline that retrieves generation code from library to produce MILP instances highly similar to target instance. MILP-Retrieval outperforms baselines in both MILP Code Generation and Instance Generation tasks, provides a novel perspective on MILP instance generation and opens new possibilities for learning-based solvers.

\end{abstract}

\begin{figure}[htbp]
    \centering
    \includegraphics[width=\linewidth]{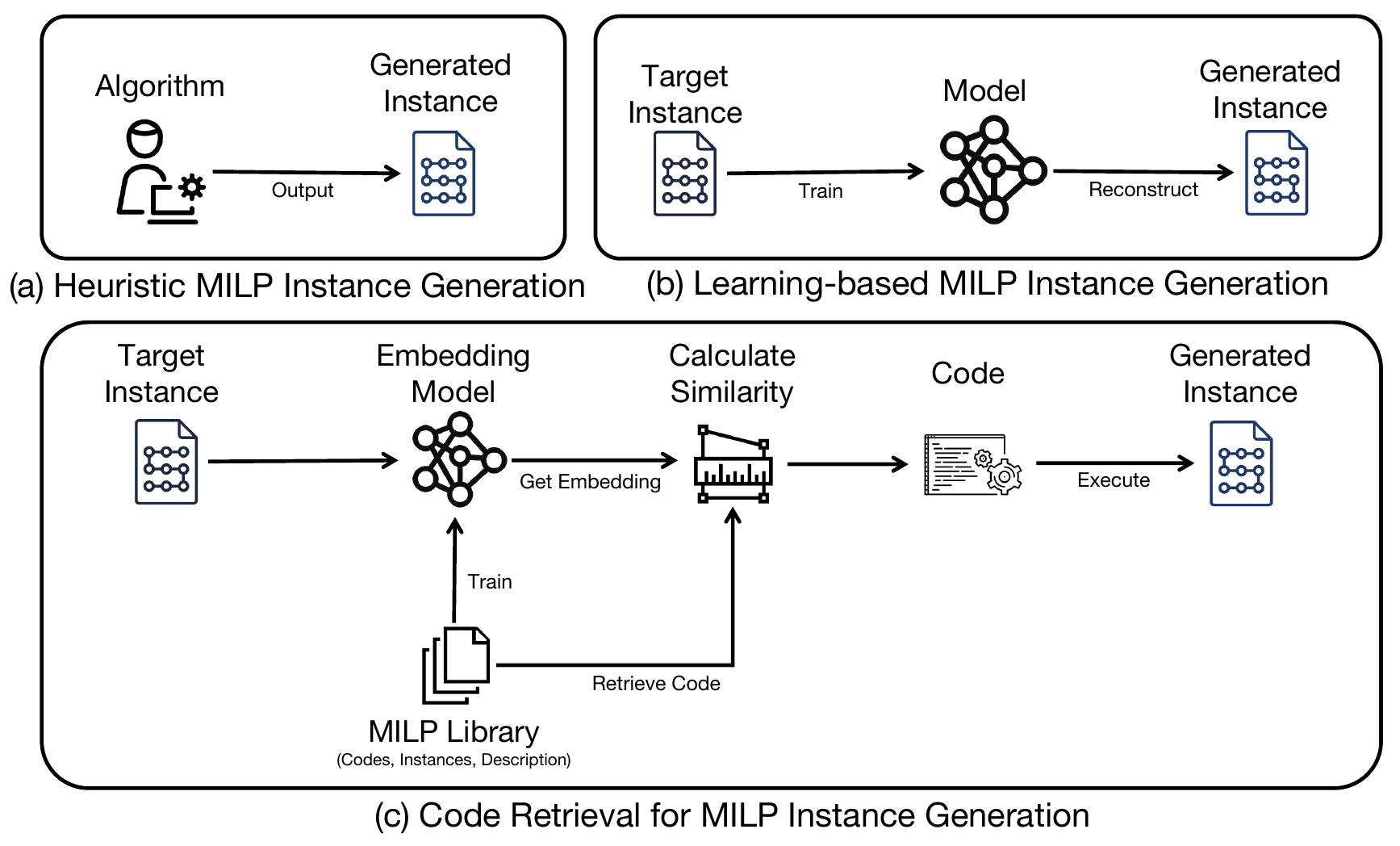}
    \caption{In MILP instance generation, (a) heuristic algorithms are used to create problem instances; (b) recent approaches train a separate model for each problem class to reconstruct problem structures; (c) our method retrieves\&tunes MILP formulation code, and uses it directly generates problem instances.}
    \label{fig:overall1}
\end{figure}

\section{Introduction}
Mixed-Integer Linear Programming (MILP) has extensive applications across various domains such as scheduling \cite{caumond2009milp, floudas2005mixed}, logistics \cite{song2018persistent, galvez2015reverse}, planning \cite{ren2010milp}, etc. Learning-based solvers \cite{li2024machine, wang2023learning, ye2023gnn} for MILP problems have shown promising performance compared to traditional solvers \cite{gurobi, BolusaniEtal2024OO, holmstrom2009user}, offering the potential for significant advancements in solving these complex problems. However, the development of learning-based solvers faces a critical challenge: data scarcity \cite{gleixner2021miplib, bengio2021machine}. Unlike domains like text or image processing, where large datasets can be scraped from the web \cite{dubey2024llama}, large-scale MILP problem datasets are not readily available. As a result, the task of MILP instance generation has been gaining increasing attention in recent years.

Early research on MILP instance generation primarily focused on using domain knowledge or heuristics to generate problem instances. For example, some approaches aimed to create problems with specific mathematical formulations and represent them using MILP \cite{rejowski2004efficient, morales2013tight, moretti2021detailed}, while others generated problems by sampling from an encoding space and using statistical data to guide the generation process\cite{smith2015generating, bowly2020generation}. These methods often relied heavily on expert knowledge to define problem templates or generation rules. This reliance on domain expertise limits their applicability in downstream tasks, such as ML-based solvers and solver hyperparameter tuning \cite{li2024machine}.

Recent research in the field has focused on learning-based approaches to MILP instance generation. These include techniques such as restructuring the structure of the MILP problem \cite{yang2024learning, liu2024milp}, generating partial structures of MILP problems based on the VAE paradigm \cite{geng2023deep, guoacm}, and using diffusion methods to reconstruct the constraints of MILP problems \cite{zhangmilp}. These approaches learn from specific classes of problems and aim to partially reconstruct problems within the same problem class. While these methods have shown great potential, they have several limitations. Each method requires retraining a model for every new problem class, which is time-intensive and computationally expensive. 

In this paper, we present a novel insight: the space of MILP problems is smaller than that of other modalities such as text, images, or videos. This smaller space allows for the creation of a sufficiently large MILP problem library. From this library, it is possible to retrieve problems that are highly similar to a given MILP problem. This insight raises three critical questions: \textbf{Q1.} \emph{How can we construct a sufficiently large MILP Problem Library?} \textbf{Q2.} \emph{How can we retrieve problems similar to a target problem from the library?} \textbf{Q3.} \emph{Why is this approach more efficient than learning-based methods for instance generation?}

To address these questions, we reformulate the original MILP Instance Generation task as \textbf{MILP Code Generation} task, which focuses on generating code capable of producing new problems similar to a given target problem. Compared to traditional MILP Instance Generation task, MILP Code Generation offers the following advantages: (1) It significantly reduces the time and computational resources needed to generate new problems; (2) It provides full control over the scale and complexity of the generated problems by allowing parameters in the code to be easily adjusted; (3) The generated problems come with corresponding mathematical formulations, ensuring interpretability.

Our approach begins by leveraging the latest method, MILP-Evolve \cite{li2024towards}, to generate a diverse MILP library, where each entry includes the MILP instance itself, its natural language description, and the code used to generate it. Next, we pretrain a MILP embedding model from the library. This model aligns the embeddings of MILP problems with their corresponding textual descriptions, drawing inspiration from the CLIP training paradigm \cite{radford2021learning}, which aligns images and text. Based on this embedding model, we introduce \textbf{MILP-EmbedSim}, a novel similarity metric for MILP problems. Unlike existing structural similarity metrics used in MILP instance generation \cite{geng2023deep, guoacm}, MILP-EmbedSim measures similarity by comparing embedding differences, enabling it to accurately capture similarities between problems of the same class but varying scales. We also propose a new pipeline, \textbf{MILP-Retrieval}, which leverages MILP-EmbedSim to identify MILP problems in the library that are similar to a given target problem. A key advantage of MILP-Retrieval is its ability to retrieve not only the problem itself but also the code used to generate it. This makes it possible to efficiently generate new problems simply by executing the retrieved code.

In the \emph{MILP Code Generation} task, we implemented two baselines due to the lack of existing methods: (1) a few-shot learning approach for code generation using GPT-4o as the large language model (LLM); (2) an LLM fine-tuned on the MILP library described earlier. MILP-Retrieval outperforms both baselines, achieving higher code validity and generating problems that more closely resemble the target problem. In the \emph{MILP Instance Generation} task, MILP-Retrieval shows comparable or superior performance across various metrics when compared to the latest open-source methods. MILP-Retrieval also significantly reduces the time and computational complexity for instance generation.

\textbf{The main contributions of the paper are as follows.}
\begin{enumerate}
    \item We reformulate the original MILP Instance Generation task as \textbf{MILP Code Generation} task, which offers a more efficient and flexible approach to generating MILP instances. The generated problems are accompanied by their corresponding mathematical formulations, ensuring interpretability.
    \item We propose \textbf{MILP-EmbedSim}, a novel similarity metric for MILP instances that accurately measures the similarity between problems of the same class but different scales, addressing limitations of previous metrics.
    \item We introduce \textbf{MILP-Retrieval}, a pipeline for constructing a MILP library and retrieving both problems and their corresponding code. MILP-Retrieval outperforms existing baselines in both the MILP Code Generation and MILP Instance Generation tasks.
\end{enumerate}

\section{Preliminary}

\subsection{MILP Problem and its Bipartite Graph Representation}

The standard formulation of a Mixed-Integer Linear Programming (MILP) problem is expressed as follows:

\begin{equation}
    \begin{aligned}
    \mathop{\min}\limits_x \quad & c^T x, \\
    \text{subject to} \quad & Ax \leq b, \\
    & l \leq x \leq u, \\
    & x_i \in \mathbb{Z}, \quad i \in \mathbb{I}.
    \end{aligned}
    \label{formula_2}
\end{equation}

Here, \( A \in \mathbb{R}^{m \times n} \) denotes the coefficient matrix that captures the structure of the MILP problem, \( b \in \mathbb{R}^m \) is the vector of constraint right-hand sides, and \( c \in \mathbb{R}^n \) represents the coefficients of the objective function. The variable bounds are specified by \( l \in (\mathbb{R} \cup \{-\infty\})^n \) and \( u \in (\mathbb{R} \cup \{+\infty\})^n \), which define the lower and upper limits of the variables, respectively. The subset of variables that must take integer values is indicated by \( \mathbb{I} \subseteq \{1, 2, \dots, n\} \).

Bipartite graph representation provides a lossless representation of MILP problems \cite{gasse2019exact}. Variables \( \mathcal{V} = \{v_1, v_2, \dots, v_n\} \) and constraints \( \mathcal{C} = \{c_1, c_2, \dots, c_m\} \) of the MILP problem are modeled as two distinct sets of nodes. An edge \( e_{ij} = (v_i, c_j) \in \mathcal{E} \) exists in the graph if the variable \( v_i \) appears in constraint \( c_j \). This results in a bipartite graph \( \mathcal{G} = (\mathcal{V}, \mathcal{C}, \mathcal{E}) \), which encodes the structural relationships between variables and constraints. For further details on the graph features utilized in this paper, please refer to Appendix \ref{bgr-detail}.

\subsection{Different Forms of MILP Data}

In addition to the Bipartite Graph Representation, this paper leverages other forms of MILP data. These forms are described below.

\textbf{Construction Code}. Construction code, referred to as code, represents the MILP problem in a generative format. In our work, this is implemented as Python scripts utilizing the pyscipopt library \cite{BolusaniEtal2024OO}. Each code snippet corresponds to a specific class of MILP problems, encapsulating the logic necessary to construct instances within that class.

\textbf{Textual Description}. Textual descriptions provide a natural language representation of MILP problems, generated using methods provided by \cite{li2024towards}. The process begins by feeding the construction code into LLM to extract key information about each class of MILP problems, such as the formulation method and associated topics. Next, the LLM is provided with statistical details about specific MILP instances. Finally, the textual description is generated, combining the broader formulation context with instance-specific statistics.

\begin{figure}[htbp]
    \centering
    \includegraphics[width=\linewidth]{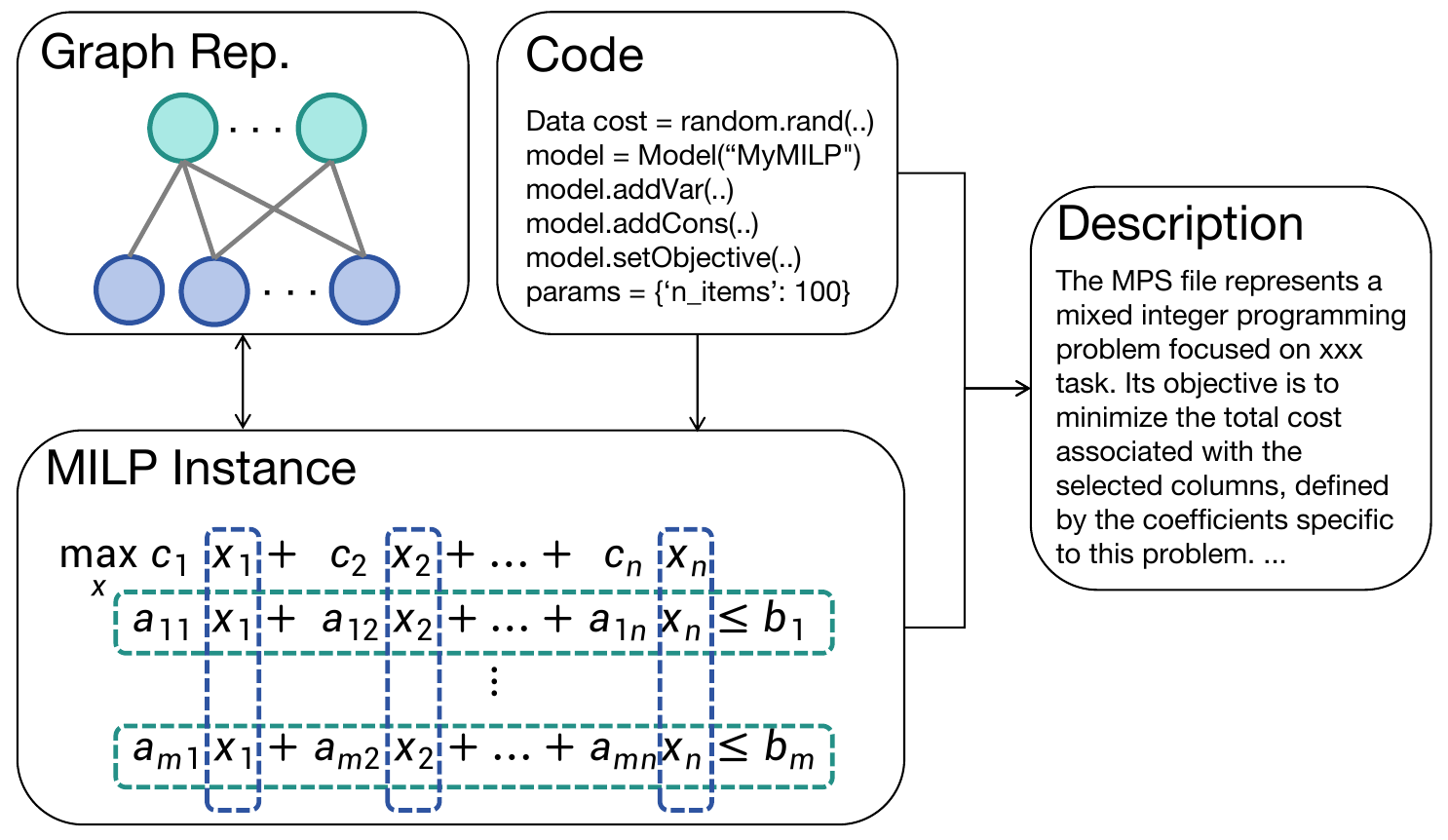}
    \caption{Relationships among the different forms of MILP data.}
    \label{fig:data-flow}
\end{figure}

Figure \ref{fig:data-flow} illustrates the relationships among the different forms of MILP data. Additionally, a sample of different forms of MILP data is provided in Appendix \ref{milpdata-sample} for reference.

\subsection{MILP Code Generation Task}

In prior work on the \emph{MILP Instance Generation} task, learning-based methods \cite{geng2023deep, yang2024learning, guoacm, zhangmilp} typically follow a paradigm where a model is trained on a dataset of MILP instances \( P = \{p_1, p_2, \dots, p_n\} \) belonging to a specific class. Once trained, the model is used to reconstruct a testing set of MILP instances \( Q = \{q_1, q_2, \dots, q_m\} \), producing a reconstructed instance set \( Q' = \{q_1', q_2', \dots, q_m'\} \). The objective of this task is to ensure that the distribution of \( P \) is as similar as possible to that of \( Q' \). For example, prior work \cite{geng2023deep} evaluates similarity by computing the Jensen-Shannon (JS) divergence \cite{lin1991divergence} between the structural statistics of instances in \( P \) and \( Q' \).

In this paper, we reformulate the \emph{MILP Instance Generation} task as the \emph{MILP Code Generation} task. Under this new formulation, a single unified model is trained on MILP problems and associated data across multiple classes, rather than being restricted to a single class. For a testing set of MILP instances \( P = \{p_1, p_2, \dots, p_n\} \), the model generates a piece of code \( c \). Executing \( c \) directly produces the instance set \( Q' = \{q_1', q_2', \dots, q_m'\} \). The objective remains the same: to minimize the divergence between the distributions of \( P \) and \( Q' \).

\begin{figure*}[htbp]
    \centering
    \includegraphics[width=\linewidth]{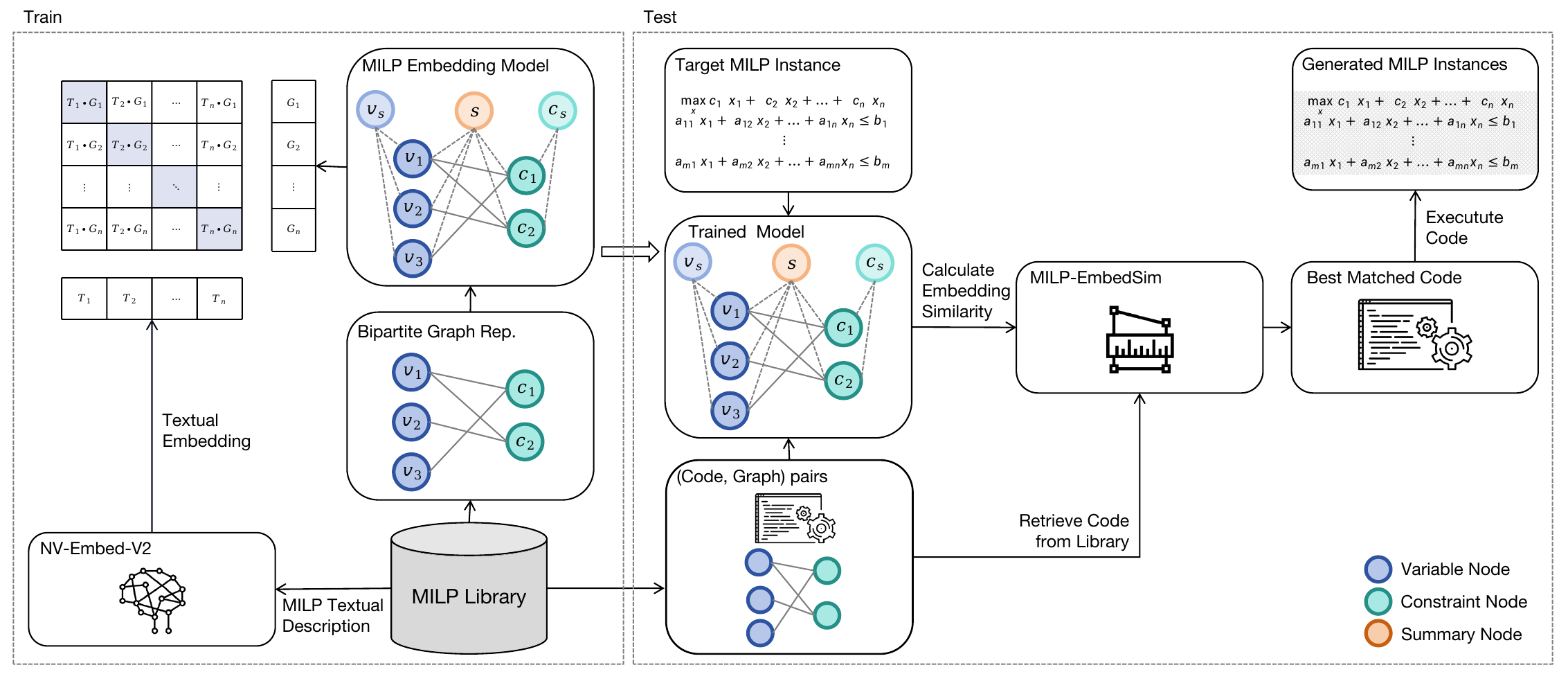}
    \caption{We begin by constructing the MILP libraries Evolve/Train and Evolve/Test. Using Evolve/Train, we train an MILP embedding model following the paradigm of contrastive learning. Based on this model, we propose a novel similarity metric called MILP-EmbedSim. Furthermore, we reformulate MILP Instance Generation as MILP Code Generation, where MILP-Retrieval retrieves code from Evolve/Train that can generate problems similar to the target problem. The retrieved code is then directly executed to generate the desired problem instances.}
    \label{fig:overall2}
\end{figure*}

This reformulated MILP Code Generation paradigm offers several significant advantages:
\begin{itemize}
    \item \textbf{Unified Model}: Unlike MILP Instance Generation, which requires training separate models for each class of MILP problems, MILP Code Generation facilitates the development of a unified model capable of handling multiple classes.  
    \item \textbf{Computational Efficiency}: By directly generating executable code to produce instances, MILP Code Generation bypasses the need for iterative structural predictions, significantly reducing computational overhead.
    \item \textbf{Full Control over Problem Characteristics}: By adjusting parameters within the generated code, users can achieve fine-grained control over the scale and difficulty of the generated problems.
\end{itemize}

This reformulation not only simplifies the MILP instance generation process but also opens new avenues for scalable and explainable MILP generation.

\section{Methodology}

In this section, we introduce the key components of our approach, which are closely interconnected. As illustrate in Figure \ref{fig:overall2}, We begin by constructing a comprehensive MILP library. This library includes a diverse collection of data types, such as MILP instances, source code, graphs, and textual description (Section \ref{sec4-1}). Leveraging the graph and textual descriptions from the MILP library, we pretrain a MILP embedding model, this model employs contrastive learning to map MILP instances into a shared embedding space, capturing their structural and semantic characteristics (Section \ref{sec4-2}). Based on the pretrained embedding model, we introduce a novel metric, MILP-EmbedSim, to quantify the similarity between MILP instances (Section \ref{sec4-3}). Leveraging this metric, we present MILP-Retrieval, a pipeline for retrieving code from the library and generating MILP instances, serves as an effective solution for the MILP Code Generation task (Section \ref{sec4-4}).

\subsection{MILP Library Preparation}
\label{sec4-1}

To construct the MILP library, we adopt method from MILP-Evolve \cite{li2024towards}. This process systematically generates a diverse set of MILP instances by evolving from initial seed classes, ensuring a rich variety of problems suitable for training and evaluation. Specifically, we prepare two separate libraries, termed Evolve/Train and Evolve/Test, dedicated to training and testing purposes, respectively. 

The Evolve/Train library is generated by evolving from 8 distinct MILP seed classes. This process yields a comprehensive dataset comprising 4,000 MILP codes and 59,033 corresponding MILP instances, graphs, and textual descriptions. This library serves as the primary dataset for training the MILP embedding model. The Evolve/Test library, on the other hand, evolves from 4 disjoint MILP seed classes, resulting in a smaller dataset of 50 MILP codes and 672 associated MILP instances, graphs, and textual descriptions. Further details on the construction process, including the specific seed classes, are provided in Appendix \ref{milplib-detail}.

\subsection{Pretraining MILP Embedding Model}
\label{sec4-2}

Using the MILP library described above, we train a powerful MILP embedding model capable of capturing both structural and semantic information. Specifically, we adopt a contrastive training framework inspired by CLIP \cite{radford2021learning}, aligning the bipartite graph representation of MILP instances with their corresponding textual descriptions. This alignment enables the model to learn a shared embedding space that effectively captures semantic relationships between different representations of MILP problems. One key advantage of this approach is that the embedding model learns to extract semantic information while ensuring that MILP problems of varying scales within the same class are mapped to similar embeddings. This consistency arises from the shared textual information associated with problems of the same class.

\begin{figure}[t]
    \centering
    \includegraphics[width=0.23\textwidth]{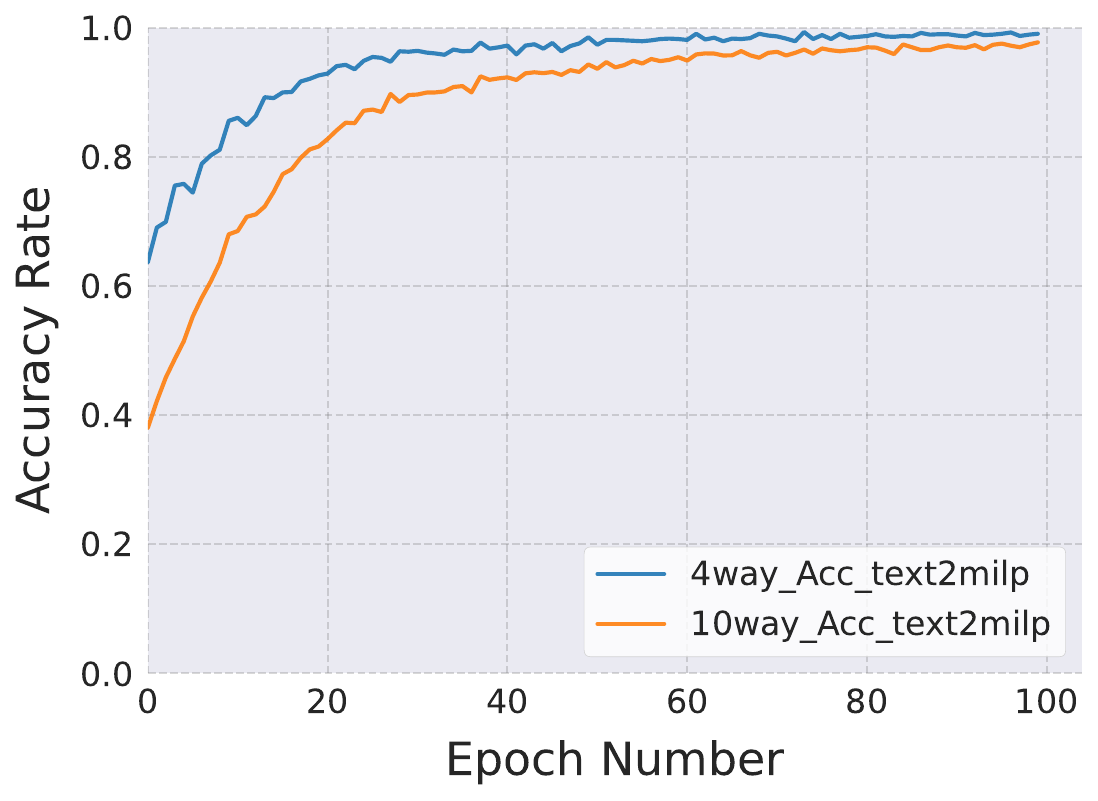}
    \hfill
    \includegraphics[width=0.23\textwidth]{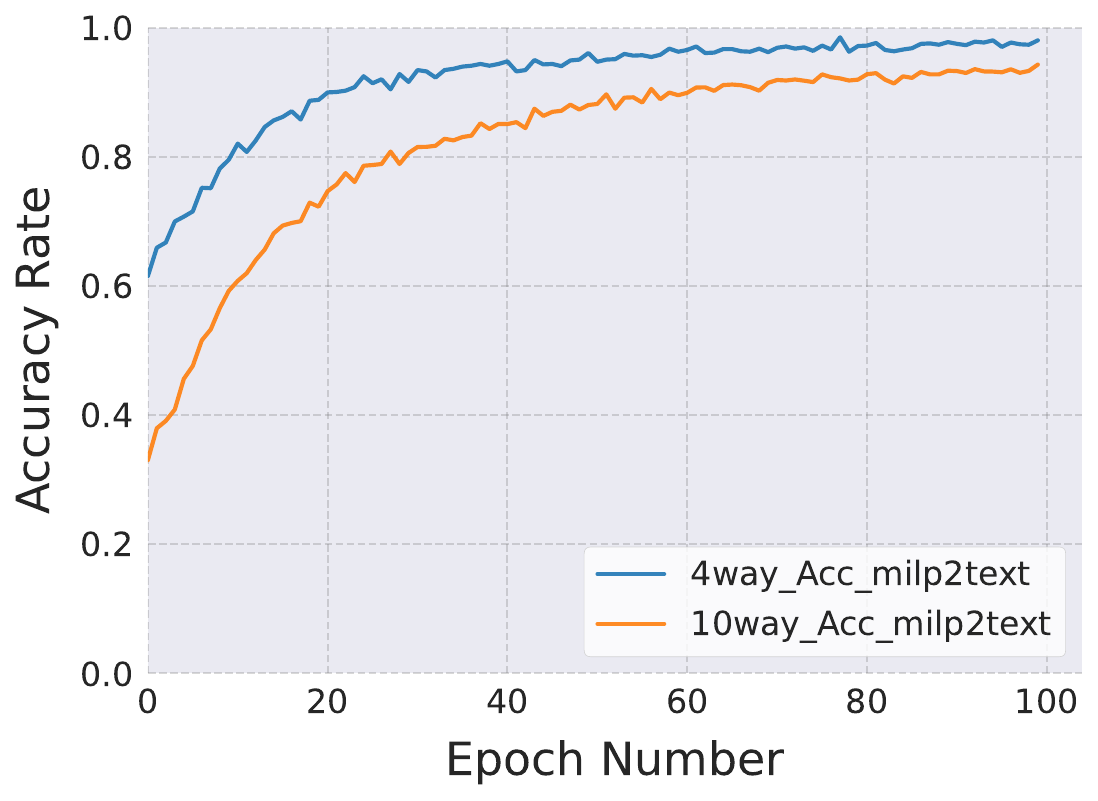}
    \caption{The Text-to-MILP and MILP-to-Text Accuracy Rate curves with respect to epoch number.}
    \label{embed-result1}
\end{figure}

Our goal is to train an MILP embedding model $f_\theta:\mathcal P\rightarrow \mathbb R^d$, where $\mathcal P$ is the space of MILP problems. For the textual embedding component $g_\theta:\mathcal T\rightarrow \mathbb R^d$, we utilize the state-of-the-art text embedding model NV-Embed-V2 \cite{lee2024nv}, freezing its weights during training. The training process employs a symmetric cross-entropy loss designed to encourage higher similarity for correct (graph, text) pairs compared to all incorrect pairings.

We divide the Evolve/Train library into training and validation sets in a 9:1 ratio. Figure \ref{embed-result1} reports the 4-way and 10-way accuracy curves on validation set during training, which measure whether the model can correctly identify the textual description corresponding to a given MILP problem (and vice versa) from 4 or 10 candidates. These results serve as intermediate results of MILP embedding model. Additional details about the MILP embedding model, including its architecture and the derivation of the loss function, are provided in the Appendix \ref{embed-detail}.

\subsection{MILP-EmbedSim: An MILP Instance Similarity Metric}
\label{sec4-3}

Existing methods for measuring similarity between groups of MILP instances often rely on computing JS divergence of their structural statistics \cite{geng2023deep, guoacm}. However, these approaches have significant limitations: they depend on a small, manually selected set of statistical indicators (only 11 metrics) and are restricted to comparing similarity between two groups of MILP problems, making pairwise comparisons infeasible.

We were inspired by the Fréchet Inception Distance (FID) \cite{heusel2017gans, salimans2016improved}, a metric used in the image generation domain to evaluate the quality of generated images, which employs Inception-V3 \cite{szegedy2016rethinking}. Based on this, we proposed MILP-EmbedSim, an MILP problem similarity metric. Experiments on the milp-to-text tasks using the MILP Embedding Model demonstrated that the model acts as an MILP classifier, with the embeddings encoding semantic information about the problem classes. MILP-EmbedSim uses the trained MILP embedding model to compute the cosine similarity between normalized embedding vectors. Formally, let $P$ and $Q$ represent two MILP instances whose similarity is to be evaluated, and $f_\theta$ denote the MILP embedding model. The MILP-EmbedSim calculation is as follows:

\begin{equation}
    x_p=\frac{f_\theta(P)}{||f_\theta(P)||},x_q=\frac{f_\theta(Q)}{||f_\theta(Q)||}
\end{equation}
\begin{equation}
    \text{MILP-EmbedSim}(P,Q)=x_px_q^T
\end{equation}

Unlike prior metrics, MILP-EmbedSim supports pairwise comparisons and offers a crucial advantage: the ability to accurately assess the similarity between instances of different scales but belonging to the same problem class. This feature arises from the design of the MILP embedding model, which is trained on data containing multiple instances of varying scales for each MILP class. These instances share similar textual descriptions, enabling the embedding model to capture cross-scale relationships.

\begin{figure*}[htbp] 
    \centering
    \begin{subfigure}[b]{0.34\textwidth} 
        \centering
        \includegraphics[width=\textwidth]{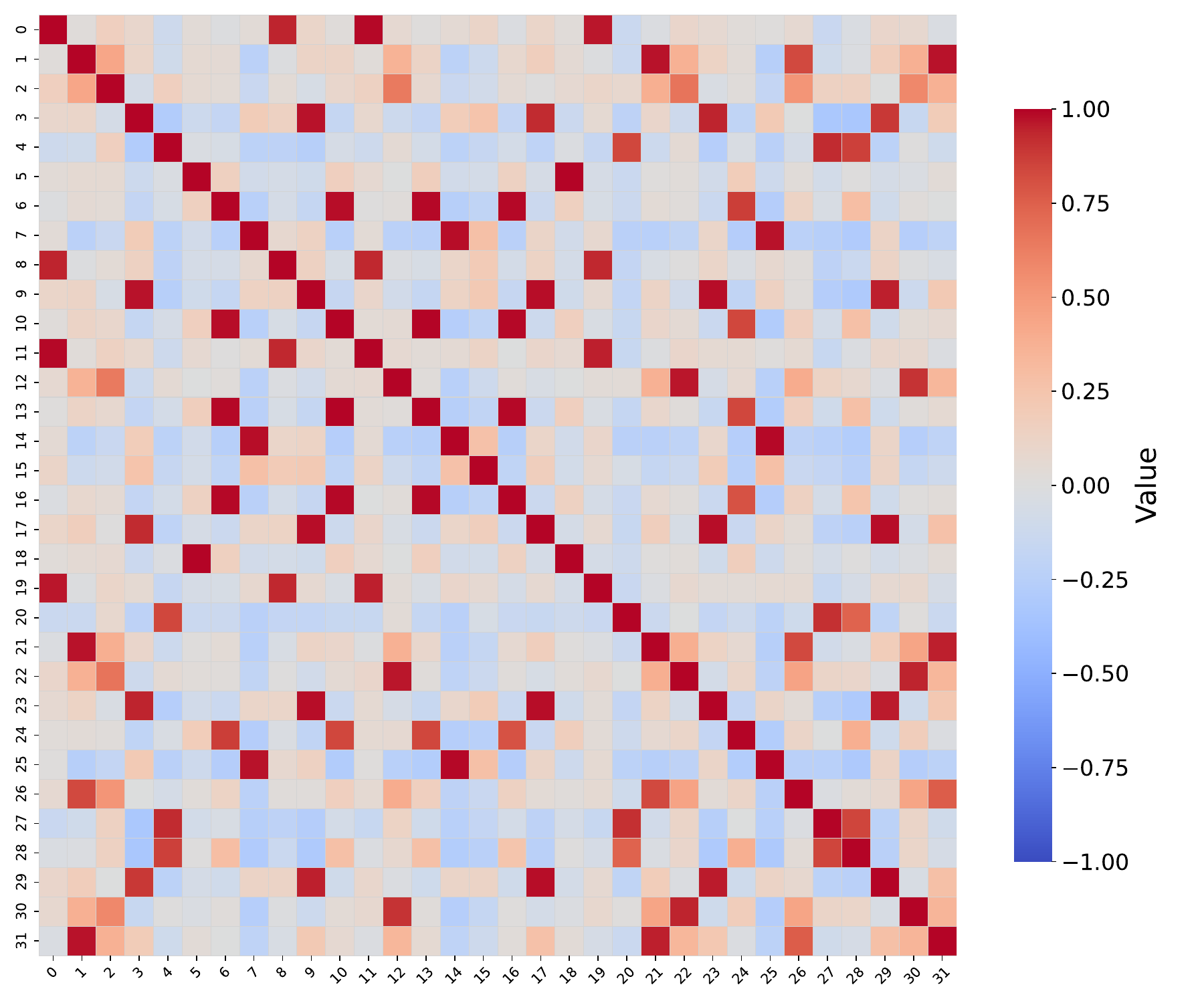} 
        \caption{Similarity between first 32 MILP classes of Evolve/Train.}
        \label{fig:vis-1}
    \end{subfigure}
    \hfill
    \begin{subfigure}[b]{0.34\textwidth}
        \centering
        \includegraphics[width=\textwidth]{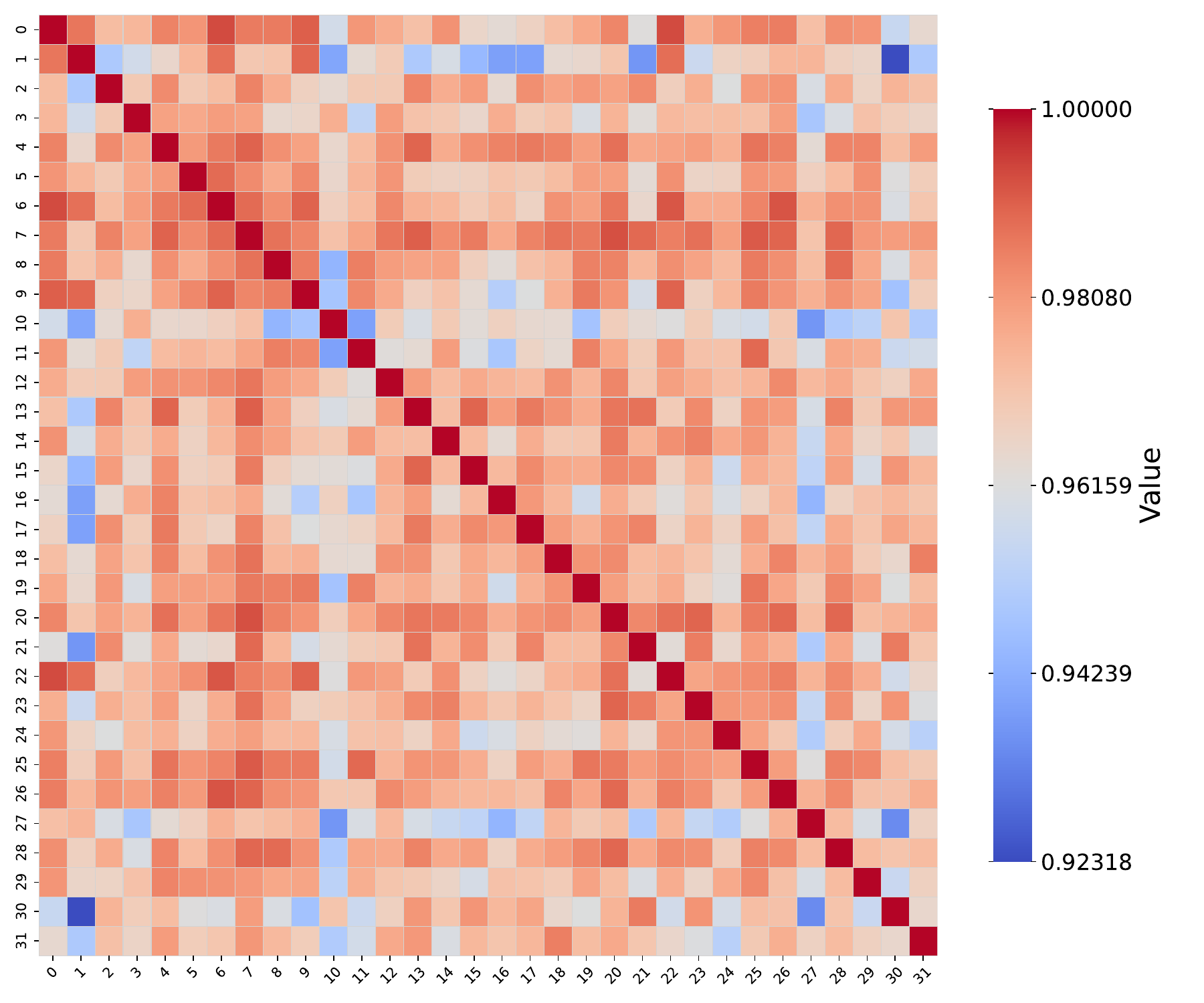}
        \caption{Similarity between 32 TSP instances of different sizes.}
        \label{fig:vis-2}
    \end{subfigure}
    \hfill
    \begin{subfigure}[b]{0.29\textwidth}
        \centering
        \includegraphics[width=\textwidth]{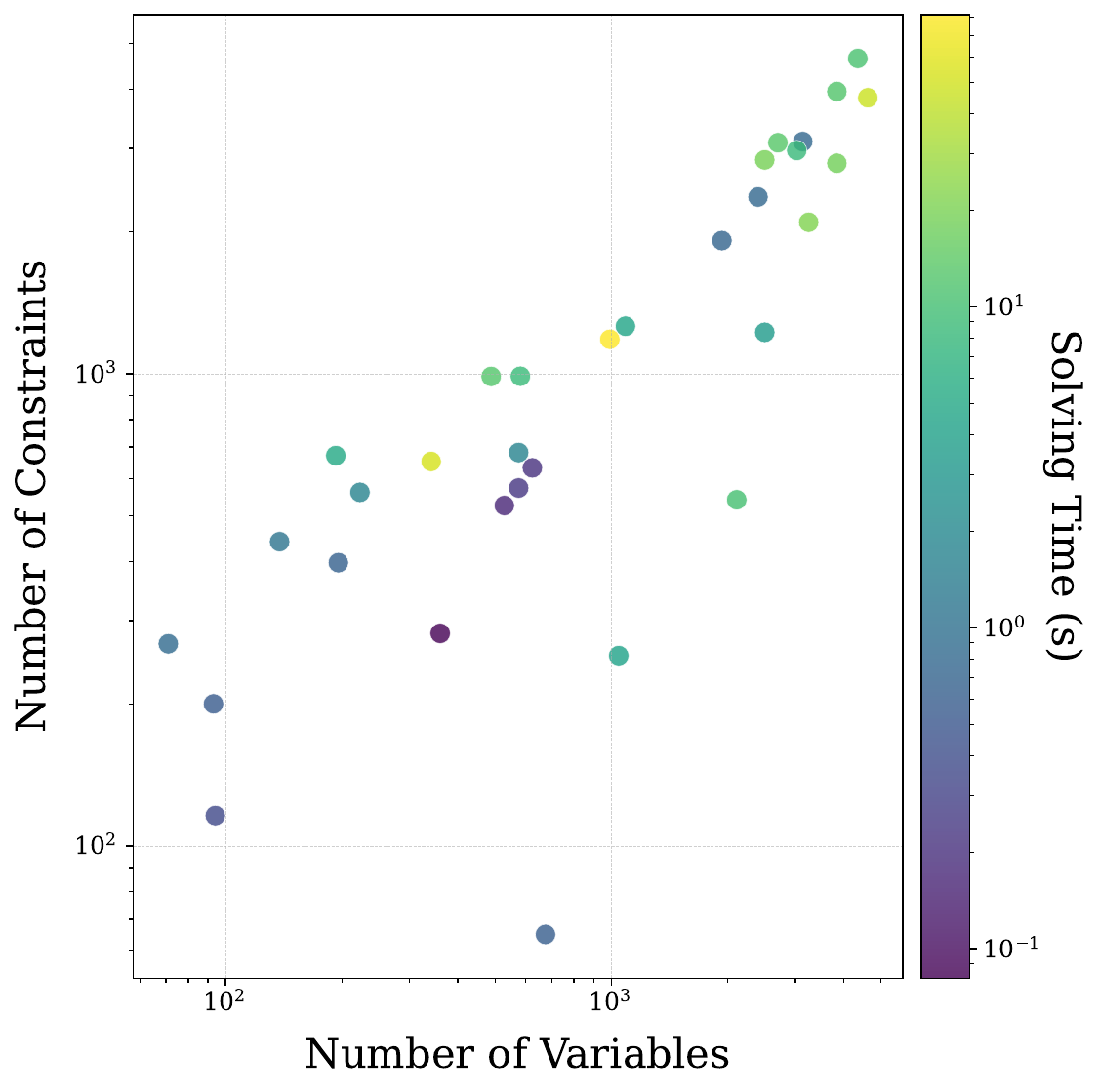}
        \caption{Problem size and solving time of 32 TSP instances.}
        \label{fig:vis-3}
    \end{subfigure}

    \caption{Visualization results of Evolve/Train dataset and TSP instances.}
    \label{fig:vis}
\end{figure*}

To illustrate the effectiveness of MILP-EmbedSim, we provide a similarity matrix in Figure \ref{fig:vis-1} for the first 32 MILP problem classes in the Evolve/Train dataset. Each matrix entry is computed using a single instance generated from the code of each class. A detailed explanation of these results is included in Appendix \ref{explain}. Additionally, we generated 32 instances of the Traveling Salesman Problem (TSP), which serves as a seed MILP class of the Evolve/Test dataset and does not appear in the training set. In Figure \ref{fig:vis-3}, we visualize the characteristics of these instances, including their scale and solving times, while Figure \ref{fig:vis-2} shows their similarity matrix. Our results demonstrate that MILP-EmbedSim generalizes effectively to unseen instances, providing robust similarity measurements even for out-of-distribution MILP problems. Further experimental results supporting these claims are provided in Appendix \ref{addition-embedsim}.

\subsection{MILP-Retrieval: Code Retrieval for MILP Instance Generation}
\label{sec4-4}

Leveraging the MILP library we constructed and the simplicity and efficiency of the MILP-EmbedSim similarity metric, we propose MILP-Retrieval, using code retrieval for MILP instance generation. Given a target MILP instance $p$, our method retrieves the most relevant code $c_k$ from MILP library $\{(p_i,c_i)\}_{i=1}^N$, where $p_i$ represents the $i$-th instance and $c_i$ represents the corresponding code for generating that instance. The retrieval process identifies $c_k$ as: $k=\text{argmax}_k \text{MILP-EmbedSim}(p,p_k)$.

Executing $c_k$ generates new MILP instances $\{q_1', q_2',\dots,q_m'\}$, effectively replicating or approximating the structural and semantic characteristics of the target instance. This design offers two significant advantages:

\begin{itemize}
    \item \textbf{Training-Free Retrieval}: Unlike previous learning-based MILP instance generation approaches, which require training a separate model for each problem class, MILP-Retrieval only relies on comparing MILP embeddings. This eliminates the need for problem-specific model training, reducing complexity and resource requirements.
    \item \textbf{Single-Instance Sufficiency}: MILP-Retrieval requires only one target instance for retrieval, in contrast to learning-based methods that require a large number of instances from the same class for training \cite{geng2023deep, guoacm}. This makes our approach highly effective in data-scarce scenarios.
\end{itemize}

\section{Experiments}

We first introduce the experimental settings in Section \ref{5-1}, including the benchmarks, baselines, and metrics. Results for the MILP Code Generation task and MILP Instance Generation task are presented in Sections \ref{5-2} and \ref{5-3}, respectively. 


\subsection{Experimental Setup }
\label{5-1}


\textbf{Benchmarks.} To ensure a fair comparison, we conducted experiments on two datasets: the Evolve/Test dataset, which contains 50 classes of problems, and a selection of instances from the MIPLIB benchmark \cite{gleixner2021miplib}. This approach accounts for the fact that some combinatorial optimization problems in the benchmark may already appear within the seed MILP classes of the Evolve/Train dataset. Both the Evolve/Train and Evolve/Test datasets are generated using MILP-Evolve \cite{li2024towards} but feature disjoint sets of MILP seed classes to avoid overlap and ensure generalization. For experiments on MIPLIB, we selected three problem categories: NurseScheduling, cvs and iis. These three classes of problems have been used in previous work on MILP instance generation \cite{geng2023deep, wang2024digmilp}.

\textbf{Baselines.} For the \emph{MILP Code Generation} task, as no existing baseline methods are available, we implemented two comparative approaches that generate MILP code from textual descriptions of problems. The first method, termed `GPT-4o', uses GPT-4o \cite{hurst2024gpt} as the underlying LLM and follows a few-shot learning paradigm. The second method, termed `Finetuned LLaMA3-8b', employs the LLaMA3-8b model \cite{dubey2024llama} fine-tuned on (textual description, code) pairs from the Evolve/Train dataset. For the \emph{MILP Instance Generation} task, we use a learning-based baseline `ACM-MILP' that adopts the Variational Autoencoder (VAE) paradigm \cite{guoacm}. This method is based on prior work \cite{geng2023deep}, modifies the original method by enabling the reconstruction of multiple constraints in a single step. We use the open source code provided by the authors of the paper to implement ACM-MILP. In our experiments, the hyperparameters for all types of problems are the same as those used for the CA (Combinatorial Auction) dataset in the original paper. Detailed implementation specifics of all baseline methods are provided in Appendix \ref{baseline-detail}. 

\begin{table*}[htbp]
\caption{Comparison between MILP-Retrieval and baseline methods on the MILP Code Generation task.}
\vskip -0.1in
\label{table: code-table1}
\begin{center}
\begin{small}
\begin{tabular}{lcccc}
\toprule
\multirow{2}{*}{Method}                & \multicolumn{3}{c}{Code Validity} & \multirow{2}{*}{Averaged Similarity} \\
                & Pass@1    & Pass@4    & Pass@10   &                                      \\
\midrule
MILP-Retrieval  & \textbf{50/50}     & -         & -         & \textbf{0.920}                                 \\
GPT-4o          & 2/50      & 6/50      & 16/50     & 0.444                                \\
Finetuned LLaMA 3-8b & 14/50     & 26/50     & 35/50     & 0.271  \\ 
\bottomrule        
\end{tabular}
\end{small}
\end{center}
\vskip -0.1in
\end{table*}

\begin{table*}[htbp]
\caption{On the MILP Instance Generation task, MILP-Retrieval is compared with the baseline method on 7 types of problems. ``Timelimit'' means that the training process cannot be completed within 72 hours, and ``Infeasible'' means that the generated problems are infeasible. The reconstruction ratio $\eta$ of ACM-MILP is set to 0.05, 0.1, and 0.2 respectively.}
\vskip -0.1in
\label{table: instance-table1}
\begin{center}
\begin{small}
\begin{tabular}{cccccc}
\toprule
\multirow{2}{*}{Problem Source} & \multirow{2}{*}{Problem Class} & \multirow{2}{*}{MILP-Retrieval} & \multicolumn{3}{c}{ACM-MILP} \\
                        &                                &                                 & $\eta=0.05$        & $\eta=0.1$     & $\eta=0.2$     \\
\midrule
\multirow{4}{*}{Evolve/Test}             & FCNF                           &  0.819                               &   0.435       &   0.397      &  0.354       \\
                        & TSP                            &   0.981                              &   \multicolumn{3}{c}{Infeasible}           \\
                        & GA                             &   0.842                              &   -0.002       &  0.015       &  -0.046       \\
                        & VRP                            &   0.992                              &  \multicolumn{3}{c}{Infeasible}           \\
\midrule
\multirow{3}{*}{MIPLIB} & NurseSched                     &   0.962                              &  0.025        &   -0.038      &  0.021       \\
                        & CVS                            &   0.822                              &     0.129     &  0.055       &  0.023       \\
                        & IIS                            &   0.919                              &     \multicolumn{3}{c}{Timelimit}         \\
\bottomrule   
\end{tabular}
\end{small}
\end{center}
\vskip -0.1in
\end{table*}

\textbf{Metrics.} We employed multiple metrics to comprehensively evaluate our proposed approach. \emph{MILP instance similarity} serves as a shared metric for both the MILP Code Generation and MILP Instance Generation tasks, which is assessed using MILP-EmbedSim. For the \emph{MILP Code Generation} task, we evaluated \emph{Code Validity}, which checks wheth er the generated code is executable, and the feasibility of the problem encoded by the generated MILP code. For the \emph{MILP Instance Generation} task, we focused on two primary metrics: \emph{computational hardness} of the generated instances and the \emph{feasible ratio}, defined as the proportion of feasible instances out of the total generated instances. 

\begin{figure}[tbp]
    \vskip -0.2in
    \centering
    \includegraphics[width=\linewidth]{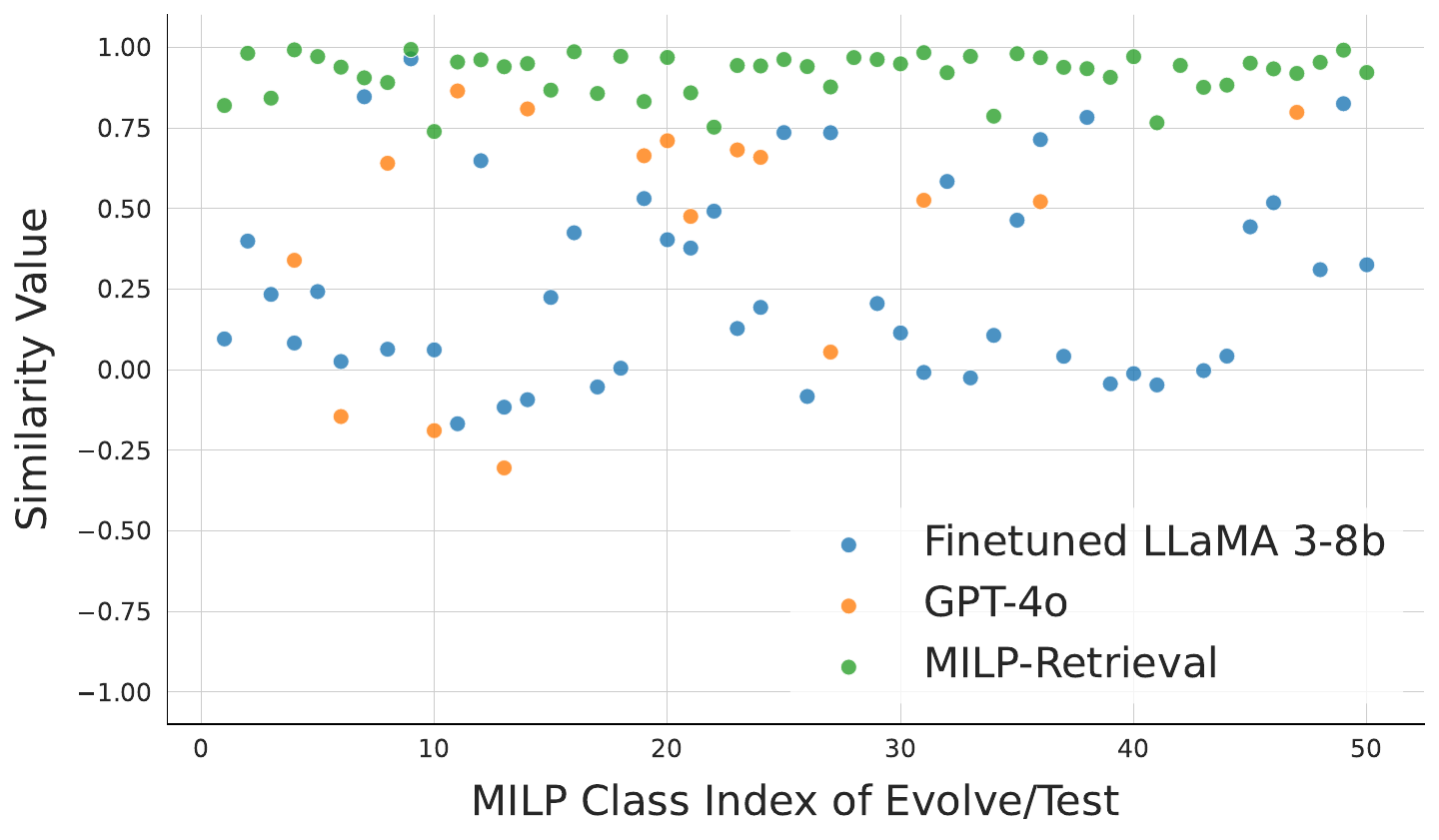}
    \caption{Results on Evolve/Test for MILP Code Generation task.}
    \label{fig:scat-plot}
    \vskip -0.1in
\end{figure}

\begin{figure}[htbp]
    \vskip -0.2in
    \centering
    \includegraphics[width=\linewidth]{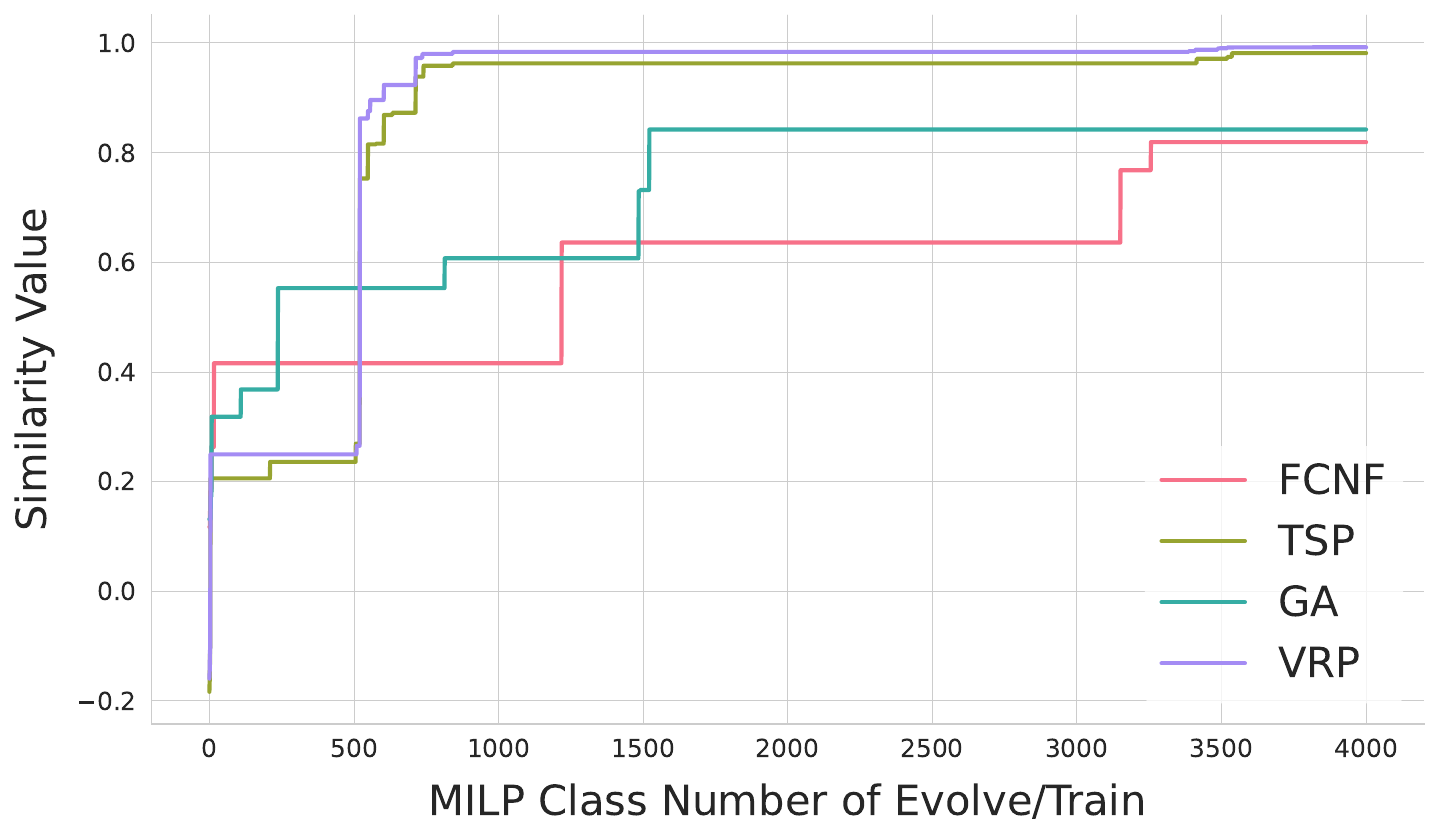}
    \caption{Ablation experiment on formulation code library size.}
    \label{fig:class-comparison}
    \vskip -0.1in
\end{figure}

\subsection{Results on MILP Code Generation Task}
\label{5-2}

Table \ref{table: code-table1} presents a comparison of our proposed MILP-Retrieval method with the two baseline approaches, evaluated on the Evolve/Test dataset comprising 50 MILP classes. Since baseline methods use textual descriptions of MILP problems as input, we did not conduct experiments on the MIPLIB benchmark. For Code Validity, we report pass@k, which measures the ratio of cases where at least one valid solution is generated after repeating the experiment $k$ times. Specifically, we conducted 10 repeated trials for both GPT-4o and Finetuned LLaMA 3-8b. For the average similarity metric, where the best results of the baseline methods within 10 trials are considered. MILP classes for which no valid code could be generated are excluded from the average similarity calculation.

\begin{figure}[htbp] 
    \centering
    \begin{subfigure}[b]{0.22\textwidth} 
        \centering
        \includegraphics[width=\textwidth]{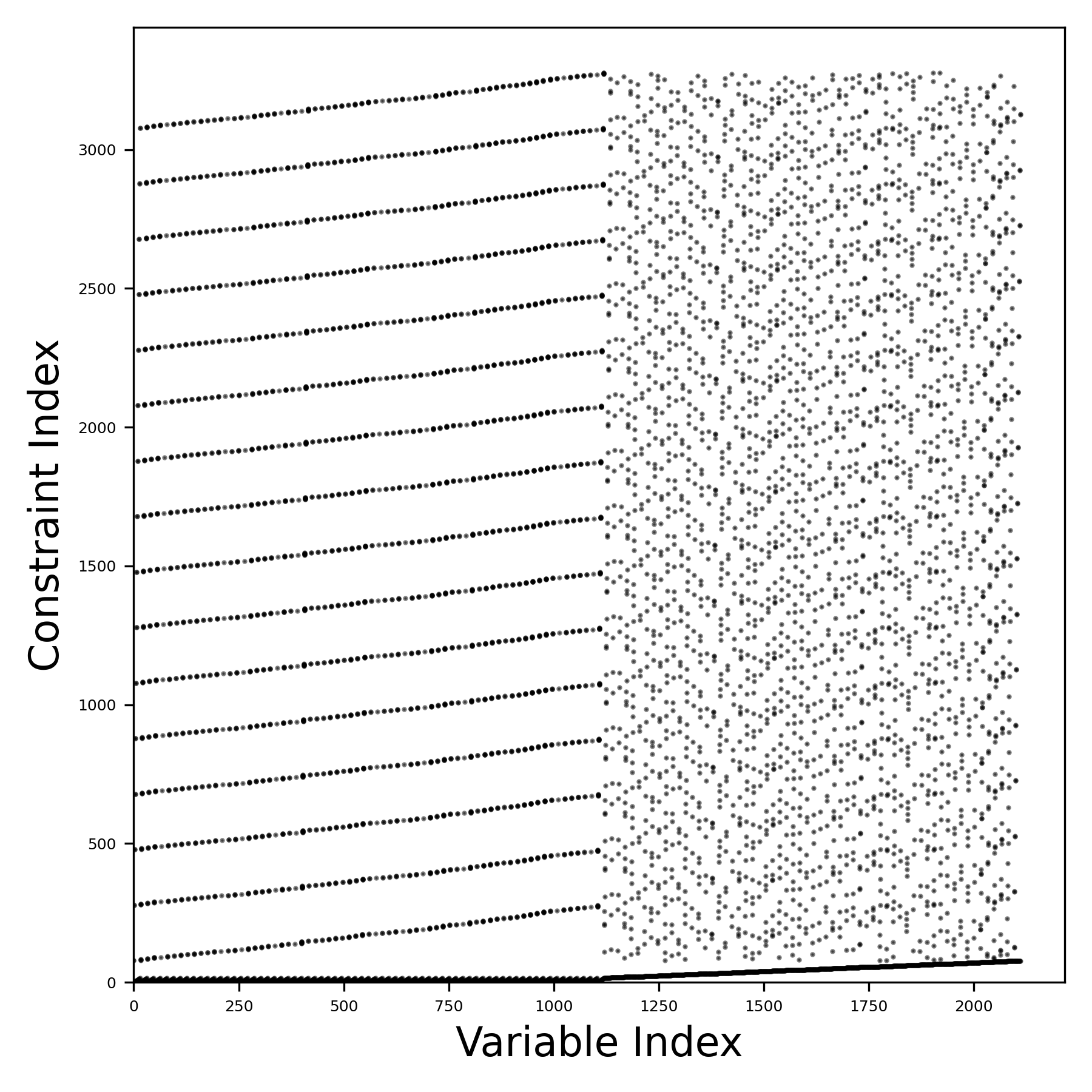} 
        \caption{Original Instance.}
        \label{fig:vis-21}
    \end{subfigure}
    \hfill
    \begin{subfigure}[b]{0.22\textwidth}
        \centering
        \includegraphics[width=\textwidth]{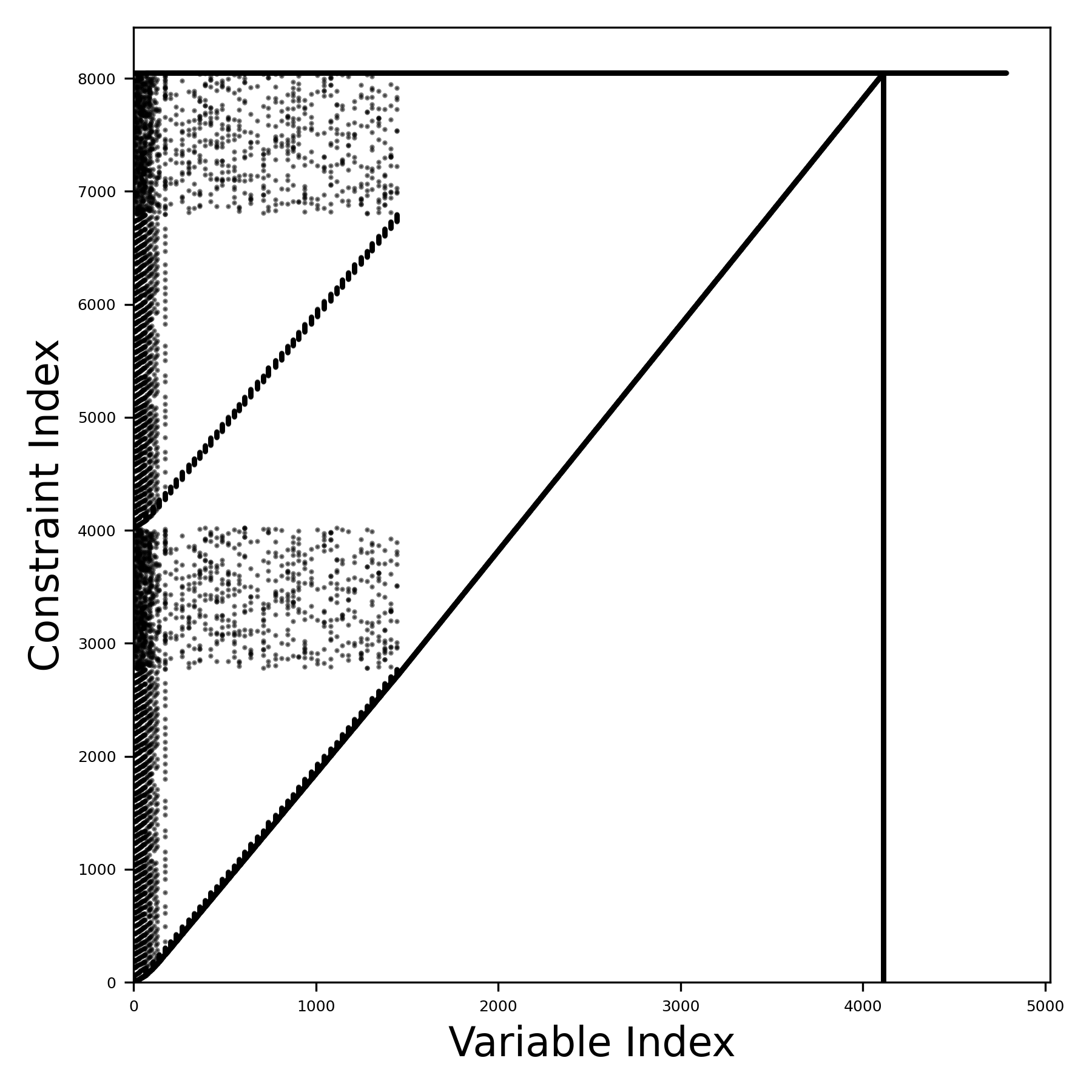}
        \caption{MILP-Retrieval.}
        \label{fig:vis-22}
    \end{subfigure}
    \hfill
    \begin{subfigure}[b]{0.22\textwidth}
        \centering
        \includegraphics[width=\textwidth]{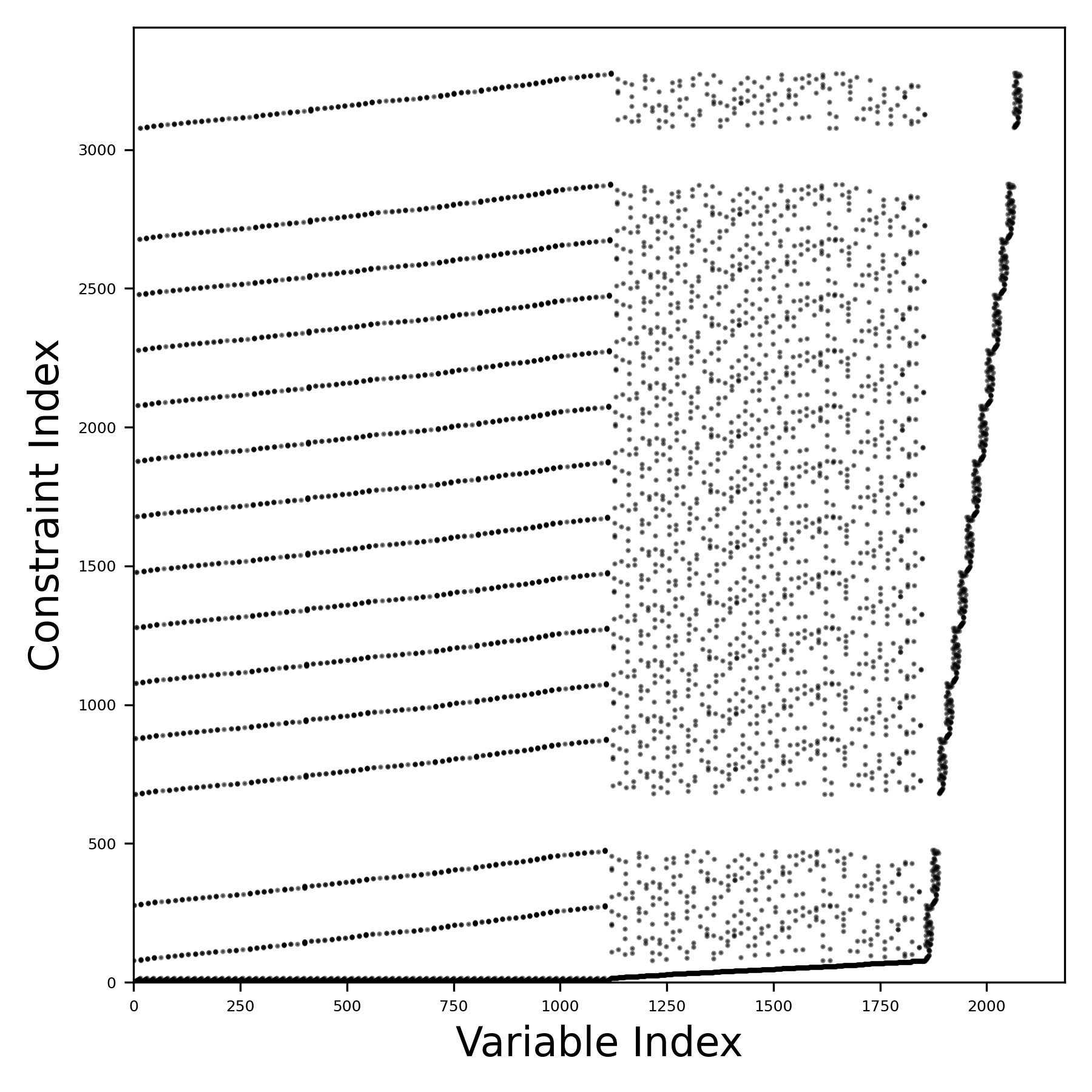}
        \caption{ACM-MILP, $\eta=0.1$.}
        \label{fig:vis-23}
    \end{subfigure}
    \begin{subfigure}[b]{0.22\textwidth}
        \centering
        \includegraphics[width=\textwidth]{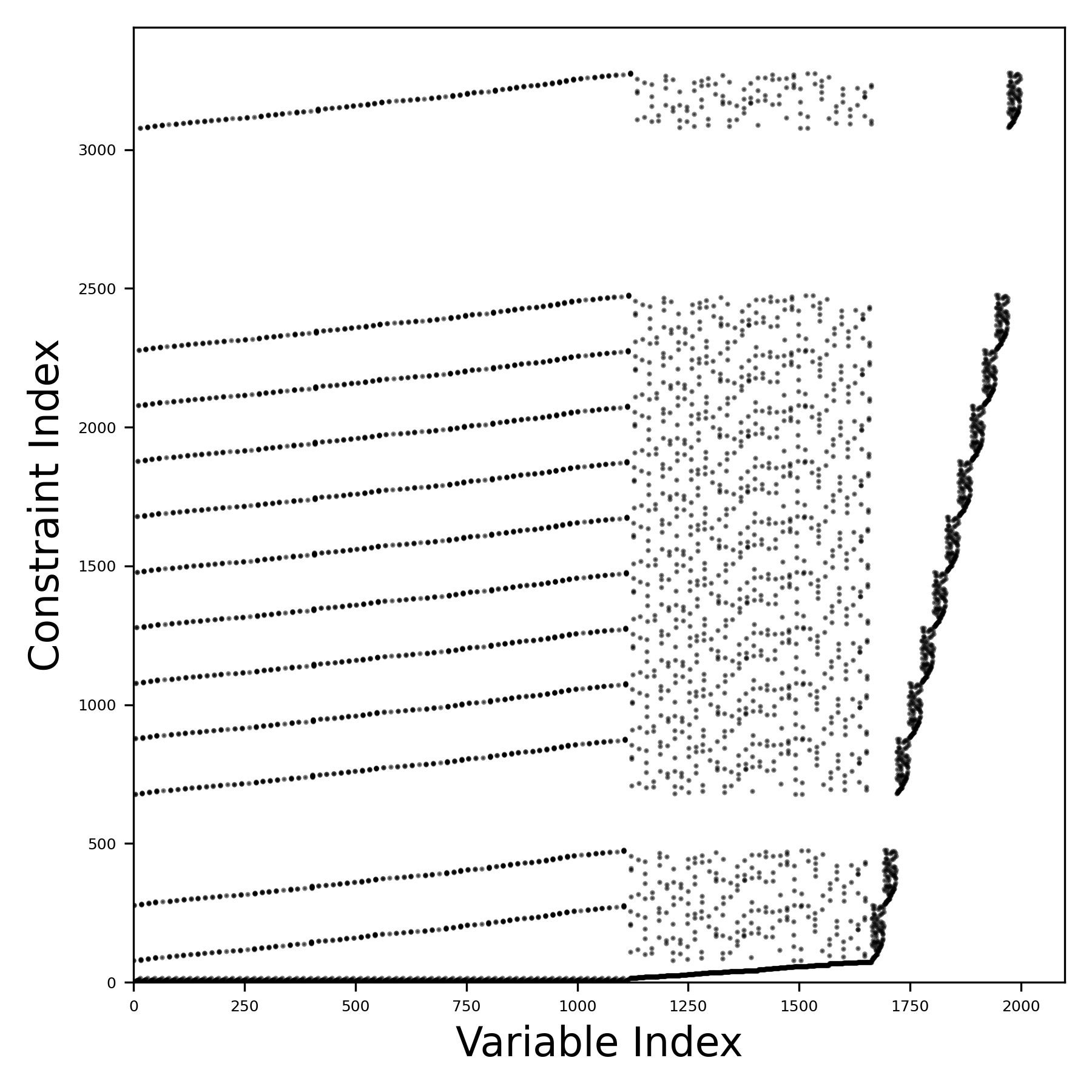}
        \caption{ACM-MILP, $\eta=0.2$.}
        \label{fig:vis-24}
    \end{subfigure}

    \caption{Visualization results of generated instances.}
    \label{fig:vis2}
\end{figure}

Figure \ref{fig:scat-plot} illustrates the best MILP instance similarity achieved by different methods across the 50 MILP classes in Evolve/Test. For GPT-4o, results are derived from 10 repeated trials, while for Finetuned LLaMA 3-8b, we report the optimal outcomes from 25 repeated experiments. Figure \ref{fig:class-comparison} highlights the retrieval performance of MILP-Retrieval on four benchmark problems (seed classes from Evolve/Test, not appearing in Evolve/Train). The results demonstrate the effectiveness of Evolve/Train, as MILP-Retrieval retrieves instances from MILP libraries of varying sizes while restricting the class number from Evolve/Train during retrieval. These experiments underscore the capability of MILP-Retrieval to leverage the library’s diversity to achieve robust and scalable results in MILP code generation.

\subsection{Results on MILP Instance Generation Task}
\label{5-3}

Given that the learning-based MILP generation method requires training a separate model for each problem class, we focus our comparison on four seed classes (FCNF, TSP, GA, VRP) from the Evolve/Test dataset and three problem classes (Nursesched, CVS, IIS) from MIPLIB. In Table \ref{table: instance-table1}, we compare the similarity of the problems generated by MILP-Retrieval and the baseline method \cite{guoacm} to the original problems. We provide dataset details of the MILP Instance Generation experiments in Appendix \ref{detail-instance-gen}.

Since MILP-Retrieval generates problems from code and allows for control over problem size and difficulty via parameter adjustments, we do not directly compare solving times. We note that the instances reconstructed by ACM-MILP exhibit patterns that differ from the original problem, resulting in lower similarity. In Figure \ref{fig:vis2}, we visualize the cvs16r70-62 problem from the CVS dataset, presenting the original instance, the instance generated by MILP-Retrieval, and the coefficient matrices of the instances generated by ACM-MILP at different reconstruction ratios.








\section{Related Works}

\subsection{Machine Learning on MILP}

Machine learning methods have demonstrated superior performance over traditional algorithms in solving various combinatorial optimization problems due to their ability to capture the characteristics of similar problems. These approaches can be broadly classified into two categories. The first category involves integrating learning-based modules into traditional solvers by replacing or augmenting key components, such as branching \cite{gasse2019exact, gupta2020hybrid, gupta2022lookback}, cut selection \cite{tang2020reinforcement, wang2023learning}, and presolve \cite{kuang2023accelerate, liu2024l2p}. The second category focuses on improving the solution search process itself. Techniques such as predict-and-optimize \cite{han2023gnn, ye2023gnn, ye2024lightmilpopt} and large neighborhood search \cite{sonnerat2021learning, huang2023searching} utilize predictive models to guide the solver toward promising regions of the solution space, thereby enhancing efficiency and solution quality. A key challenge in both categories is the availability of sufficient MILP data for training these models. This challenge highlights the critical need for generating diverse and high-quality MILP instances.

\subsection{MILP Instance Generation}

The field of MILP instance generation has traditionally relied on heuristic methods to create problem instances tailored to specific types or statistical characteristics \cite{smith2015generating, bowly2020generation}. While effective in controlled scenarios, these methods often lack the flexibility to address broader applications or more diverse instance distributions. Learning-based MILP generation methods use model to  learn the distribution of the problems and reconstruct them. For example, some methods focus on restructuring the problem’s underlying structure \cite{liu2024milp, yang2024learning}, while others utilize paradigms like VAE or diffusion models to reconstruct problem constraints \cite{geng2023deep, wang2024digmilp, guoacm, zhangmilp}. Recent work \cite{li2024towards} proposes a novel approach for generating diverse MILP problems. Our approach MILP-Retrieval, along with the concept of MILP Code Generation, offers a novel perspective on MILP instance generation.

\section{Conclusion}

In this paper, we introduced a novel approach, using code retrieval for MILP instance generation. We began by training an MILP embedding model, which served as the foundation for MILP-EmbedSim, a similarity metric designed to quantify the resemblance between MILP instances. Leveraging this metric, we reformulated the MILP Instance Generation task as an MILP Code Generation problem, leading to MILP-Retrieval—a retrieval-based framework that identifies and executes code capable of generating new problem instances similar to a given target problem.




\section*{Impact Statement}

Although this paper is related to generative AI, it primarily focuses on the field of combinatorial optimization, discussing a new paradigm for MILP instance generation. Compared to works on generating images, texts, videos, and other content, our current work does not involve sensitive aspects such as privacy, bias, or factual accuracy. Therefore, we believe there are no societal consequences that need to be specifically highlighted in this context.


\bibliographystyle{icml2025}

\newpage
\appendix
\onecolumn

\section{Details of Bipartite Graph Features}
\label{bgr-detail}

To encode an MILP instance as a corresponding bipartite graph, we incorporate information about both variables and constraints into the node features of the graph representation. The specific node features used in our encoding are detailed in Table \ref{tab:node_features}.

Additionally, the bipartite graph features include solution-related information about the MILP instance. To obtain this data, we solve each problem instance using Gurobi \cite{gurobi}, with a computation time limit of 50 seconds per instance. This ensures a standardized and practical approach to extracting solution-based features while maintaining computational efficiency.

\begin{table}[htbp]
\caption{Node type features and descriptions for Variables and Constraints.}
\label{tab:node_features}
\begin{center}
\begin{small}
\begin{tabular}{ccc}
\toprule
\textbf{Node Type} & \textbf{Feature} & \textbf{Description} \\
\midrule
        \multirow{9}{*}{\textbf{Vars}} & norm coef & Objective coefficient, normalized by objective norm \\
                     & type & Var type (binary, integer, impl. integer, continuous), one-hot \\
                     & has lb & Lower bound indicator \\
                     & has ub & Upper bound indicator \\
                     & solval & Solution value \\
                     & solfrac & Solution value fractionality \\
                     & sol\_is\_at\_lb & Solution value equals lower bound \\
                     & sol\_is\_at\_ub & Solution value equals upper bound \\
                     & basestat & Simplex basis status (lower, basic, upper, zero), one-hot \\
\midrule   
        \multirow{7}{*}{\textbf{Cons}} & rank & Rank of a row \\
                     & norm\_nnzrs & Fraction of nonzero entries \\
                     & bias & Unshifted side normalized by row norm \\
                     & row\_is\_at\_lhs & Row value equals left hand side \\
                     & row\_is\_at\_rhs & Row value equals right hand side \\
                     & dualsol & Dual LP solution of a row, normalized by row and objective norm \\
                     & norm\_intcols & Fraction of integral columns in the row \\
\bottomrule
\end{tabular}
\end{small}
\end{center}
\end{table}

\section{Sample of Different Forms of MILP Data}
\label{milpdata-sample}

Here we provide a sample of code and textual description in MILP Data, which comes from the Set Cover problem and is one of the seed classes of Evolve/Train. The parameter part in the code can be adjusted to control the size and solving time of the generated problem.

\textbf{Code}:

\begin{lstlisting}[language=Python]
import random
import time
import scipy
import numpy as np
import networkx as nx
from pyscipopt import Model, quicksum

class SetCover:
    def __init__(self, parameters, seed=None):
        for key, value in parameters.items():
            setattr(self, key, value)

        self.seed = seed
        if self.seed:
            random.seed(seed)
            np.random.seed(seed)

    ################# Data Generation #################
    def generate_instance(self):
        nnzrs = int(self.n_rows * self.n_cols * self.density)

        # compute number of rows per column
        indices = np.random.choice(self.n_cols, size=nnzrs)  # random column indexes
        indices[:2 * self.n_cols] = np.repeat(np.arange(self.n_cols), 2)  # force at leats 2 rows per col
        _, col_nrows = np.unique(indices, return_counts=True)

        # for each column, sample random rows
        indices[:self.n_rows] = np.random.permutation(self.n_rows) # force at least 1 column per row
        i = 0
        indptr = [0]
        for n in col_nrows:
            # empty column, fill with random rows
            if i >= self.n_rows:
                indices[i:i+n] = np.random.choice(self.n_rows, size=n, replace=False)

            # partially filled column, complete with random rows among remaining ones
            elif i + n > self.n_rows:
                remaining_rows = np.setdiff1d(np.arange(self.n_rows), indices[i:self.n_rows], assume_unique=True)
                indices[self.n_rows:i+n] = np.random.choice(remaining_rows, size=i+n-self.n_rows, replace=False)

            i += n
            indptr.append(i)

        # objective coefficients
        c = np.random.randint(self.max_coef, size=self.n_cols) + 1

        # sparce CSC to sparse CSR matrix
        A = scipy.sparse.csc_matrix(
                (np.ones(len(indices), dtype=int), indices, indptr),
                shape=(self.n_rows, self.n_cols)).tocsr()
        indices_csr = A.indices
        indptr_csr = A.indptr

        res =  {'c': c, 
                'indptr_csr': indptr_csr, 
                'indices_csr': indices_csr}

        return res

    ################# PySCIPOpt Modeling #################
    def solve(self, instance):
        c = instance['c']
        indptr_csr = instance['indptr_csr']
        indices_csr = instance['indices_csr']

        model = Model("SetCover")
        var_names = {}

        # Create variables and set objective
        for j in range(self.n_cols):
            var_names[j] = model.addVar(vtype="B", name=f"x_{j}", obj=c[j])

        # Add constraints to ensure each row is covered
        for row in range(self.n_rows):
            cols = indices_csr[indptr_csr[row]:indptr_csr[row + 1]]
            model.addCons(quicksum(var_names[j] for j in cols) >= 1, f"c_{row}")

        # Set objective: Minimize total cost
        objective_expr = quicksum(var_names[j] * c[j] for j in range(self.n_cols))

        model.setObjective(objective_expr, "minimize")
        
        start_time = time.time()
        model.optimize()
        end_time = time.time()

        return model.getStatus(), end_time - start_time

if __name__ == '__main__':
    seed = 42
    parameters = {
        'n_rows': 750,
        'n_cols': 1500,
        'density': 0.05,
        'max_coef': 100,
    }

    set_cover_problem = SetCover(parameters, seed=seed)
    instance = set_cover_problem.generate_instance()
    solve_status, solve_time = set_cover_problem.solve(instance)

    print(f"Solve Status: {solve_status}")
    print(f"Solve Time: {solve_time:.2f} seconds")
\end{lstlisting}

\textbf{Textual Description}:
\begin{quote}
    The MPS file, named `SetCover', represents a mixed integer programming problem focused on a Set Cover optimization task. Its objective is to minimize the total cost associated with the selected columns, defined by the coefficients specific to this problem. The formulation leverages inequalities to ensure that each of the 750 constraints guarantees that every row is covered by at least one selected column. The decision variables are binary, reflecting the choice of each column's inclusion in the cover. The file employs a structured approach for encoding the problem, facilitating efficient solving by optimization algorithms.
\end{quote}

\section{Details of MILP Library}
\label{milplib-detail}

To facilitate the training and evaluation of our approach, we constructed two MILP libraries: Evolve/Train and Evolve/Test. The Evolve/Train library is specifically designed for training purposes, while Evolve/Test is intended for evaluation. For Evolve/Train, we selected eight problem categories as seed classes: IS, SC, CA, CFL, KS, GIS, NF, and SAT. For Evolve/Test, we designated four distinct problem categories as seed classes: FCNF, TSP, GA, and VRP. The seed classes of Evolve/Train and Evolve/Test are disjoint.

The Evolve/Train library comprises a total of 4,000 MILP classes, while Evolve/Test contains 50 MILP classes. Each MILP class generates 20 MILP instances, followed by a filtering step that removes instances that fail to yield a feasible solution within 50 seconds. After this filtering process, the final dataset includes 59,033 MILP instances in Evolve/Train and 672 MILP instances in Evolve/Test.

The MILP class evolution process was conducted using GPT-4o-mini as the LLM. Beginning with the seed classes, the evolution procedure to construct the Evolve/Train and Evolve/Test libraries took approximately four weeks and incurred a total cost of around \$50. For detailed information on the MILP class generation process, please refer to \cite{li2024towards}. The sources of the seed classes are summarized in Table \ref{tab:seed-class}.

\begin{table}[htbp]
\caption{8 Seed Classes for Evolve/Train and 4 Seed Classes for Evolve/Test.}
\label{tab:seed-class}
\begin{center}
\begin{small}
\begin{tabular}{cccc}
\toprule
\textbf{Dataset} & \textbf{Abbreviation} & \textbf{Full Name} & \textbf{Reference} \\
\midrule
        \multirow{8}{*}{\textbf{Evolve/Train}} & IS & Maximum Independent Set & \cite{bergman2016decision} \\
                     & SC & Set Cover & \cite{balas1980set} \\
                     & CA & Combinatorial Auction & \cite{leyton2000towards} \\
                     & CFL & Capacitated Facility Location & \cite{cornuejols1991comparison} \\
                     & Knapsack & Multiple Knapsack & \cite{pisinger1999exact} \\
                     & GIS & Generalized Independent Set & \cite{colombi2017generalized} \\
                     & NF & Multicommodity Network Flow & \cite{hewitt2010combining} \\
                     & SAT & Max Satisfiability & \cite{bejar2009generating} \\
\midrule   
        \multirow{4}{*}{\textbf{Evolve/Test}} & FCNF & Fixed-Charge Network Flow & \cite{kim1999solution} \\
                     & TSP & Traveling Salesman Problem & \cite{matai2010traveling} \\
                     & GA & Generalized Assignment & \cite{cattrysse1992survey} \\
                     & VRP & Vehicle Routing Problem & \cite{braekers2016vehicle} \\
\bottomrule
\end{tabular}
\end{small}
\end{center}
\end{table}

\section{Implementation Details of MILP Embedding Model}
\label{embed-detail}


\subsection{Model Architecture}

Our MILP embedding model is composed of two primary components: (1) a bipartite Graph Neural Network (GNN) designed to capture relational information between constraints (rows) and variables (columns), and (2) a Transformer-based self-attention module that refines the extracted features.

We represent each MILP instance as a bipartite graph $(\mathcal V,\mathcal C,\mathcal E)$, where $\mathcal V$ denotes nodes corresponding to variables, $\mathcal C$ denotes nodes representing constraints, and $\mathcal E$ consists of edges connecting variables to the constraints in which they appear. Following the bipartite graph feature described in Appendix \label{tab:node_features}, the node features for variables and constraints are represented as $\mathcal V\in\mathbb R^{n\times 16}$ and $\mathcal C\in\mathbb R^{m\times 7}$, respectively.

To embed the nodes into a shared latent space of dimension $\text{emb\_size}$, we employ two separate Multi-Layer Perceptrons (MLPs) for variables and constraints. Let $x_u,x_u\in\mathbb R^{\text{emb\_size}}$ represent the constraint and variable embeddings, respectively. Additionally, edge features—such as the coefficients of variables in constraints—are encoded using a linear projection layer to enrich the representation of graph connectivity.

We leverage a bipartite GNN to propagate information between constraints and variables. Each message-passing iteration consists of two steps: (1) a column-to-row update, where variable embeddings are aggregated into constraint nodes, and (2) a row-to-column update, where constraint embeddings are propagated back to variable nodes. Additionally, we introduce a summary node that connects to all other nodes, capturing a global representation of the graph.

For message passing, we utilize a Graph Convolution Module \cite{kipf2017semisupervised} as the update function. Formally, let $\mathbf{x}{u}^{(k)}$ and $\mathbf{x}{v}^{(k)}$ be the constraint and variable embeddings at layer $k$. The updates are performed as follows:

\begin{equation}
    \mathbf{x}{u}^{(k+1)} \;=\; \mathbf{x}{u}^{(k)} + \text{BipartiteConv}\Bigl(\mathbf{x}{v}^{(k)}, \mathcal E_{uv}\Bigr)
\end{equation}
\begin{equation}
    \mathbf{x}{v}^{(k+1)} \;=\; \mathbf{x}{v}^{(k)} + \text{BipartiteConv}\Bigl(\mathbf{x}{u}^{(k+1)}, \mathcal E_{uv}\Bigr)
\end{equation}

Following the GNN layers, we sample a subset of 512 node embeddings from the output. These sampled embeddings, along with the averaged embeddings of variable and constraint nodes and the summary node embedding, form a set of 515 embeddings. This set is concatenated and fed into a stack of Transformer encoder layers. The Transformer layers apply self-attention mechanisms to model long-range dependencies and enhance the global structure of the MILP instance. Finally, the output of the Transformer module is a fixed-size MILP embedding vector in $\mathbb R^{4096}$.

\subsection{Derivation of Loss Function}

Let $(\mathcal P_i, \mathcal T_i)$ for $i=1,\dots,N$ be a batch of $N$ matched MILP–text pairs, $f_\theta$ be the MILP embedding model, producing an MILP embedding $\mathbf{p}_i = f_\theta(\mathcal P_i) \in \mathbb{R}^d$, $g_\theta$ be the text encoder, producing a text embedding $\mathbf{t}_i = g_\theta(\mathcal T_i) \in \mathbb{R}^d$. Both $\mathbf{p}_i$ and $\mathbf{t}_i$ are typically L2-normalized to have unit length, $\|\mathbf{p}_i\|_2 = 1,\|\mathbf{t}_i\|_2 = 1$. For each MILP–text pair $(i, j)$ in the batch, we define the similarity score as the dot product: $s_{ij} \;=\; \,\mathbf{p}_i^\top \mathbf{t}_j$.

Our training objective is a bidirectional contrastive objective: it treats each MILP instance $\mathbf{p}_i$ as a query and tries to classify the correct text $\mathbf{t}_i$ among all texts $\{\mathbf{t}_j\}$, and symmetrically, each text $\mathbf{t}_i$ tries to classify the correct MILP instances $\mathbf{v}_i$ among all instances $\{\mathbf{v}_j\}$.

For a fixed MILP embedding $\mathbf{p}_i$, the MILP-to-text cross-entropy loss is:
\begin{equation}
    \ell_i^\text{(MILP-to-Text)}=-\log(\frac{\exp(s_{ii})}{\sum\limits_{j=1}^N\exp(s_{ij})})
\end{equation}

Similarly, for a fixed text embedding $\mathbf{t}_i$, the text-to-MILP cross-entropy loss is:
\begin{equation}
    \ell_i^\text{(Text-to-MILP)}=-\log(\frac{\exp(s_{ii})}{\sum\limits_{j=1}^N\exp(s_{ji})})
\end{equation}

To incorporate both MILP-to-text and text-to-MILP objectives, the final symmetric loss sums these two cross-entropy terms for each pair and then averages over the batch:
\begin{equation}
    \mathcal{L} \;=\; \frac{1}{2N} \sum_{i=1}^N \ell_i^\text{(MILP-to-Text)} + \ell_i^\text{(Text-to-MILP)}
\end{equation}

\subsection{Prompt Details of NV-Embed-V2}

We use the text embedding model NV-Embed-V2 \cite{lee2024nv} as $g_\theta$ and freeze its weights during training. This model is an instruction embedding model. We set the following instruction prompt:
\begin{quote}
    Instruct: Given a linguistic description, retrieve the corresponding Mixed-Integer Linear Programming problem.
\end{quote}

\subsection{Training Details}

We trained the MILP embedding model on the Evolve/Train dataset, which contains a total of 59,033 (MILP, textual description) pairs. We randomly divided it into a training set and a validation set in a ratio of 9:1, using the training set as training data. In Figure \label{embed-result1}, we report the experimental results on the validation set. The training process was completed on a single Nvidia H800 and took about 40 hours. We provide the hyperparameters used for training in Table \ref{tab:hyper}.

\begin{table}[htbp]
\centering
\caption{Hyperparameter of MILP embedding model.}
\label{tab:hyper}
\begin{tabular}{cc|cc}
\toprule
\textbf{Name} & \textbf{Value} & \textbf{Name} & \textbf{Value} \\
\midrule
Embed\_Size & 64 & Num. of GCN Layers & 2 \\
Num. of Sampled Nodes & 512 & Num. of Attention Layers & 6 \\
Embedding Space & $\mathbb R^{4096}$ & Epoch Number & 100 \\
Learning Rate & 0.001 & Batch Size & 64 \\
Num. of Attention Heads & 8 & Optimizer & Adam \\

\bottomrule
\end{tabular}
\end{table}

\section{Implementation Details of Baselines}
\label{baseline-detail}

\subsection{`GPT-4o' Baseline}

Our implemented `GPT-4o' baseline takes as input the textual description of an MILP problem and outputs code. We employ a few-shot learning approach to prompt the LLM to generate structured MILP code. Specifically, we randomly select three (textual description, code) pairs from the Evolve/Test dataset as examples, and then use the textual description from Evolve/Test as the test input. The prompt used in this process is as follows:

\begin{quote}
Please generate Python code for the Mixed-Integer Linear Programming problem corresponding to the description below.

\emph{\{target desc}\}

Sample description 1:

\emph{\{sample\_desc1}\}

Sample code 1:

\emph{\{sample\_code1}\}

Sample description 2:

\emph{\{sample\_desc2}\}

Sample code 2:

\emph{\{sample\_code2}\}

Sample description 3:

\emph{\{sample\_desc3}\}

Sample code 3:

\emph{\{sample\_code3}\}
\end{quote}

\subsection{`Finetuned LLaMA 3-8b' Baseline}

We implemented another baseline `Finetuned LLaMA 3-8b', which also takes the textual description of the MILP problem as input and outputs the code, and tests it on Evolve/Test. We use all the data from Evolve/Train to construct the SFT dataset, where each piece of data corresponds to a (textual description, code) pair. The format of the data is as follows:
\begin{lstlisting}[language=Python]
messages = [
    {"role": "system", "content": "You are an expert in Mixed-Integer Linear Programming."},
    {"role": "user", "content": "Please generate Python code for the Mixed-Integer Linear Programming problem corresponding to the description below. \n" + description},
    {"role": "assistant", "content": code}
]
\end{lstlisting}

During testing, we use the same user prompt as input and feed the code generated by the finetuned model into GPT-4o for checking, ensuring that the output code contains no syntax errors. The prompt used for the code checking is as follows:
\begin{quote}
    Identify and fix the errors in this code, then output the complete corrected code.
\end{quote}

We use the XTuner framework \cite{2023xtuner} to perform full parameter fine-tuning on LLaMA-3-8b-instruct \cite{dubey2024llama}, with the fine-tuning hyperparameters specified in Table \ref{tab:hyper2}.

\begin{table}[htbp]
\centering
\caption{Hyperparameter of Finetuning LLaMA-3-8b.}
\label{tab:hyper2}
\begin{tabular}{cc|cc}
\toprule
\textbf{Name} & \textbf{Value} & \textbf{Name} & \textbf{Value} \\
\midrule
Epoch Num & 8 & Learning Rate & 2e-5 \\
Batch Size & 1 & Accumulate Counts & 16 \\

\bottomrule
\end{tabular}
\end{table}

\section{Dataset of MILP Instance Generation Experiments}
\label{detail-instance-gen}

In the MILP Instance Generation experiment, we used two types of datasets. The first type consists of the first four MILP classes from Evolve/Test (FCNF, TSP, GA, VRP). For each problem, we generated 20 instances to be used as training data for ACM-MILP, and these 20 instances were also used as the problems reconstructed by ACM-MILP. For the three datasets from MIPLIB (Nursesched, CVS, IIS), we used all available instances provided by MIPLIB as both training and test data for ACM-MILP. The dataset statistics are provided in Table \ref{table: instance-table2}.

\begin{table*}[htbp]
\caption{Dataset Statistics of MILP Instance Generation Experiment.}
\vskip -0.1in
\label{table: instance-table2}
\begin{center}
\begin{small}
\begin{tabular}{cccccc}
\toprule
Problem Source & Problem Class & Instance Num. & Average $|\mathcal V|$ & Average $|\mathcal C|$ & Average $|\mathcal E|$\\

\midrule
\multirow{4}{*}{Evolve/Test}             & FCNF                           &  20                               &   1096       &   594      &  2192       \\
                        & TSP                            &   20                              &   1604 & 1567 & 7592           \\
                        & GA                             &   20                              &   125000 & 750 & 250000       \\
                        & VRP                            &   20                              &  1088 & 1153 & 7168          \\
\midrule
\multirow{3}{*}{MIPLIB} & NurseSched                     &   5                              &  19501        &   7231      &  373018       \\
                        & CVS                            &   5                              &   2536     &  3397       &  9150      \\
                        & IIS                            &   2                              &   256 & 7551 & 99552      \\
\bottomrule   
\end{tabular}
\end{small}
\end{center}
\vskip -0.1in
\end{table*}

\section{Further Experiments and Analysis}
\label{add-exp}

\subsection{MILP-EmbedSim}
\label{addition-embedsim}

\begin{figure}[htbp]
    \centering
    \includegraphics[width=0.5\linewidth]{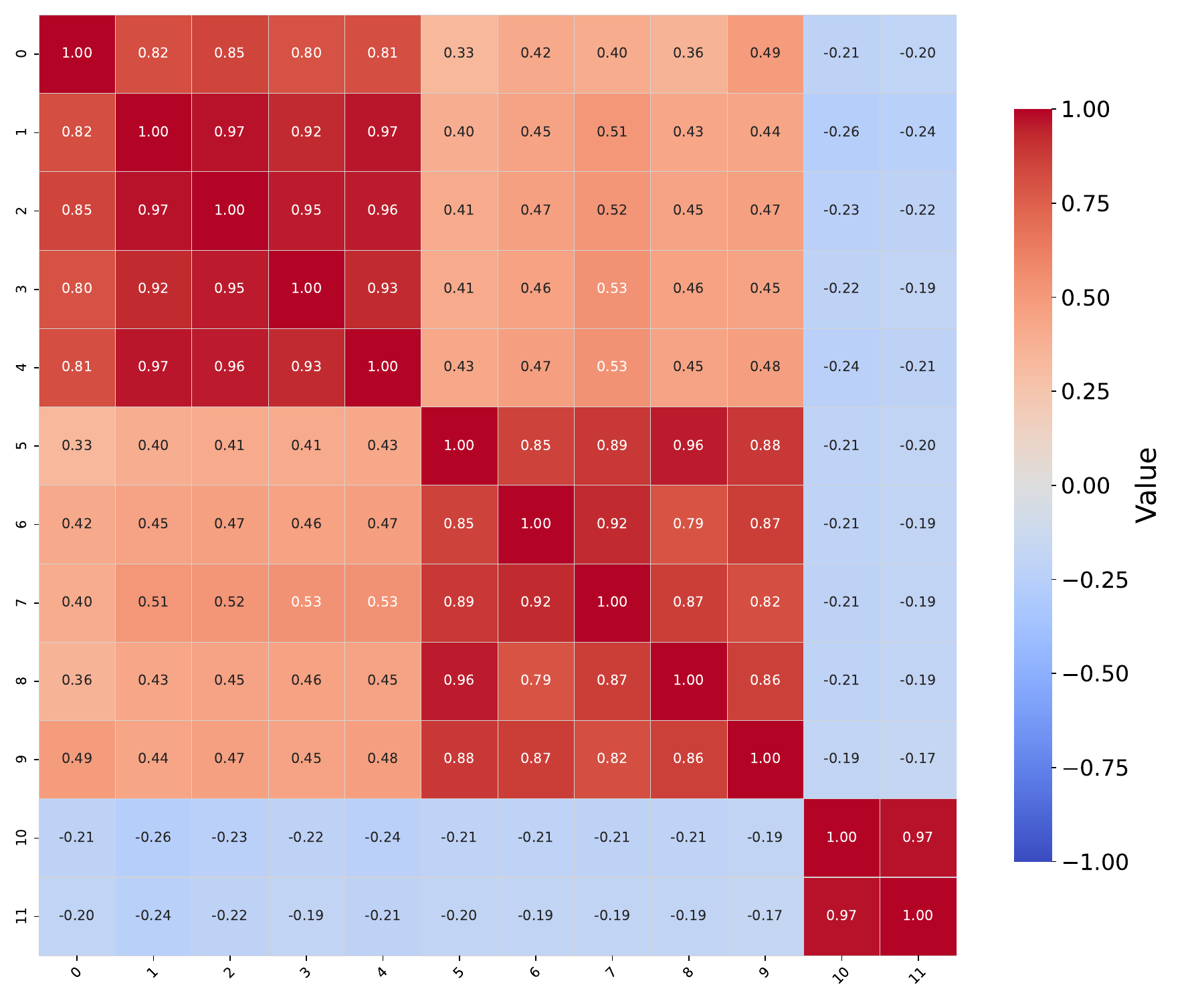}
    \caption{Similarity between instances of 3 classes from MIPLIB.}
    \label{fig:sim33}
\end{figure}

We also used MILP-EmbedSim to calculate the similarity between three types of problems from MIPLIB, further validating our similarity metric. The results are shown in Figure \ref{fig:sim33}, where the first five instances correspond to five CVS instances, the middle five instances are from the nursesched problem, and the last two instances are from the IIS problem.

\subsection{Explanation of Figure \ref{fig:vis-1}}
\label{explain}

In Figure \ref{fig:vis-1}, we visualize the similarity between the first 32 problem classes from Evolve/Train. The first 8 classes serve as seed classes, where each class is distinct from the others, resulting in relatively low similarity between them. The remaining 24 classes are evolved sequentially, and thus, there is a certain degree of similarity both among them and between these evolved classes and the seed classes. For example, problem C is a combination of problems A and B, so C retains a certain level of similarity with both A and B. This visualization also partially reveals the evolutionary path of the MILP classes.



\end{document}